\documentclass{amsart}
\usepackage{graphicx}
\usepackage{amssymb}
\usepackage{amsfonts}
\swapnumbers
\newtheorem{theorem}{Theorem}[section]

\newtheorem{corollary}[theorem]{Corollary}

\theoremstyle{definition}

\newtheorem{assumption}[theorem]{Assumption}
\newtheorem{remark}[theorem]{Remark}

\numberwithin{equation}{section}
 \theoremstyle{plain}
 
 \numberwithin{equation}{section} 
 \numberwithin{figure}{section} 
 \theoremstyle{plain}
 \theoremstyle{remark}
 \newtheorem*{acknowledgement*}{Acknowledgement}

\newcommand{\cB}{{\mathcal B}}

\newcommand{\cF}{{\mathcal F}}
\newcommand{\cG}{{\mathcal G}}
\newcommand{\cH}{{\mathcal H}}

\newcommand{\cP}{{\mathcal P}}

\newcommand{\cW}{{\mathcal W}}

\newcommand{\te}{{\theta}}
\newcommand{\vt}{{\vartheta}}
\newcommand{\Om}{{\Omega}}
\newcommand{\om}{{\omega}}
\newcommand{\ve}{{\varepsilon}}
\newcommand{\del}{{\delta}}
\newcommand{\Del}{{\Delta}}
\newcommand{\gam}{{\gamma}}
\newcommand{\Gam}{{\Gamma}}
\newcommand{\vf}{{\varphi}}

\newcommand{\Sig}{{\Sigma}}
\newcommand{\sig}{{\sigma}}
\newcommand{\al}{{\alpha}}
\newcommand{\be}{{\beta}}

\newcommand{\la}{{\lambda}}

\newcommand{\vs}{{\varsigma}}

\newcommand{\bbR}{{\mathbb R}}
\newcommand{\bbS}{{\mathbb S}}

\newcommand{\bbW}{{\mathbb W}}

\newcommand{\bbI}{{\mathbb I}}

\begin{document}
\title[]{Almost sure diffusion approximation in averaging:
direct proofs with rough paths flavors}%
 \vskip 0.1cm
 \author{ Yuri Kifer\\
\vskip 0.1cm
 Institute  of Mathematics\\
Hebrew University\\
Jerusalem, Israel}%
\address{
Institute of Mathematics, The Hebrew University, Jerusalem 91904, Israel}
\email{ kifer@math.huji.ac.il}%

\thanks{ }
\dedicatory{  }
 \date{\today}
\begin{abstract}\noindent
We consider again the fast-slow motions setups in the continuous time
$\frac {dX_N(t)}{dt}=N^{1/2} \sig(X_N(t))(\xi(tN))+b(X_N(t)),\, t\in [0,T]$ and the discrete time
$X_N((n+1)/N)=X_N(n/N)+N^{-1/2}\sig(X_N(n/N))\xi(n)+N^{-1}b(X_N(n/N)),\, n=0,1,...,[TN]$
where $\sig$ and $b$ are smooth matrix and vector functions, respectively, $\xi$ is a centered vector stationary
stochastic process with weak dependence in time and $N$ is a big parameter. We obtain estimates for the almost
 sure approximations of the process $X_N$ by certain diffusion process $\Sig$. In \cite{FK} and in other recent
 papers concerning similar setups the results were obtained relying fully on the rough paths theory. Here we derive
 our probabilistic results as corollaries of quite general deterministic estimates which are obtained with all details
 provided following somewhat ideology of the rough paths theory but not relying on this theory per se
 which should allow a more general readership to follow complete arguments.
\end{abstract}
\maketitle
\markboth{Yu.Kifer}{Almost sure diffusion approximation}
\renewcommand{\theequation}{\arabic{section}.\arabic{equation}}
\pagenumbering{arabic}

\section{Introduction}\label{sec1}\setcounter{equation}{0}

The motivation for this paper comes from the study of the asymptotic behavior as $N\to\infty$ of solutions $X_N$ of systems of ordinary
differential equations having the form
\begin{equation}\label{1.1}
\frac {dX_N(t)}{dt}=\sqrt N \sig(X_N(t))\xi(tN))+b(X_N(t)),\,\, t\in[0,T]
\end{equation}
and the corresponding discrete time problem about the asymptotic behavior of processes $X_N$ given by the recurrence relation
\begin{equation}\label{1.2}
X_N((n+1)/N)=X_N(n/N)+N^{-1/2}\sig(X_N(n/N))\xi(n)+N^{-1}b(X_N(n/N)),
\end{equation}
$n=0,1,...[TN]-1$. In both setups $\xi$ is a mean zero vector stationary process with certain weak dependence properties.
These setups include the slow and the fast motions $X_N$ and $\xi$, respectively, and the problems arising here concern the
question what happens with the slow motion as $N\to\infty$.

This kind of problem was studied already about sixty years ago in \cite{Kha} followed by \cite{PK} and \cite{Bor} 
where weak convergence of $X_N$ to a diffusion solving the stochastic differential equation of the form
\begin{equation}\label{1.3}
d\Xi(t)=\sig(\Xi(t))dW(t)+\big(b(\Xi(t))+c(\Xi(t))\big)dt
\end{equation}
was established under conditions on $\xi$ which did not allow applications to important situations such as the case when the
process $\xi$ is generated by a (deterministic) dynamical system. More recently this problem attracted a renewed interest (see
\cite{CFKMZ} and references there). These newer type of results derived above mentioned weak convergence under more general conditions
which allowed applications to dynamical systems but the proof there relied crucially on the rough paths theory. One of motivations of
the present paper is to provide arguments which can be followed by a general probability audience and not only by readers familiar with
the rough paths theory.

In \cite{FK} we extended weak convergence results to almost sure (a.s.) approximation of the process $X_N$ by an appropriately chosen
version of the diffusion $\Xi$. The main step of the proof in \cite{FK} consisted of establishing strong (i.e. a.s.) invariance principles
for the integrals $S_N(t)=N^{-1/2}\int_0^{tN}\xi(s)ds$ and sums $S_N(t)=N^{-1/2}\sum_{0\leq n<[tN]}\xi(n)$, as well as for iterated integrals
$\bbS_N(t)=N^{-1}\int_0^{tN}ds\int_0^s\xi(u)\otimes\xi(s)du$ and iterated sums $\bbS_N(t)=N^{-1}\sum_{0\leq n<m<[tN]}\xi(n)\otimes\xi(m)$.
Observe that recently such strong invariance principles were extended to all multiple iterated integrals and sums in \cite{Ki23+}. After these
was obtained we relied in \cite{FK} completely on results from the rough paths theory in \cite{FZ} concerning estimates of the Lipschitz dependence of
solutions of rough differential equations on the driving processes. This result from \cite{FZ} is proved for general processes with jumps and it
relies on notations and the machinery of the rough paths theory so that the arguments there are not easy to follow by a general probability reader.
Essentially, both in our paper \cite{FK} and in recent papers yielding only weak convergence such as \cite{CFKMZ}, the rough paths theory came as a
"black box" so that only experts in this theory could enter it to follow the complete argument. The problem we are dealing with is not really connected
with the rough paths theory, and so rather to rely on specific results from this theory, it is more natural to provide direct proofs for this problem which can be accessible to the general probability audience.

In this paper we rely on strong invariance principles for integrals, sums, iterated integrals and iterated sums proved in \cite{Ki23+} under more
general than in \cite{FK} conditions and provide a direct proof of a.s. diffusion approximations in our averaging setups which will enable a more
general readership to follow detailed arguments. In fact, the main part of the paper is devoted to a quite general deterministic setup which may have
an independent interest and from which our probabilistic results will be derived as corollaries relying mainly on \cite{Ki23+}. We start with general
vector functions $X(s,t)=X(t)-X(s)$ and $\Xi(s,t)=\Xi(t)-\Xi(s)$ without specifying their nature and work under a basic approximation assumption which
an expert in the rough paths theory would recognize as a definition of rough integrals or something similar to the, so called, sewing lemma.
It turns out that already in such general setup
it is possible to obtain appropriate estimates for $X-\Xi$ in a H\" older norm which will be applied to our specific probabilistic setup for required a.s.
approximation of $X_N$ given by (\ref{1.1}) or (\ref{1.2}) by $\Xi$ given by (\ref{1.3}). This strategy makes the proof direct and easier to follow.

In the deterministic part we will proceed somewhat similarly to \cite{FH} but observe that the book \cite{FH} does not provide estimates of the local
Lipschitz constant of the, so called, It\^ o--Lyons map which is crucial for our purposes. Moreover, the arguments for the H\" older norms employed
 in \cite{FH} cannot be applied directly to the process $X_N$ given by (\ref{1.2}) since it has jumps at the times $k/N$, but we still employ certain
  modified H\" older seminorms until the last steps, which is one of novelties of this paper, passing to variational norms at the very end. This
  allows a more straightforward proof than working with variational norms from the beginning which involves, so called, control functions leading to a rather complicated machinery requiring arguments which are
   difficult to follow by non experts in rough paths theory. We observe also that some
    estimates for discrete rough paths were obtained in \cite{BFT} but they are not appropriate for our purposes.

  The structure of the present paper is as follows. In Section \ref{sec2} we describe our main results. Sections \ref{sec3} and \ref{sec4}
  are devoted to the deterministic setups in the continuous and discrete setups, respectively, and in Section \ref{sec5} we derive our results
  for the original probabilistic setup as corollaries from the previous sections.

\section{Preliminaries and main results}\label{sec2}\setcounter{equation}{0}

\subsection{General deterministic setup: continuous time}\label{subsec2.1}

We will start with several maps $X:[0,T]\to\bbR^d$, $\Xi:[0,T]\to\bbR^d$, $S:[0,T]\to\bbR^D$, $W:[0,T]\to\bbR^D$, $\bbS:\{ 0\leq s\leq t\leq T\}\to
\bbR^D\otimes\bbR^D$,
$\bbW:\{ 0\leq s\leq t\leq T\}\to\bbR^D\otimes\bbR^D$, $\sig:\bbR^d\to\bbR^d\otimes\bbR^D$ and $b:\bbR^d\to\bbR^d$ where $D\geq d\geq 1$. It will be
convenient to write $X(s,t)=X(t)-X(s)$, $\Xi(s,t)=\Xi(t)-\Xi(s)$, $S(s,t)=S(t)-S(s)$ and $W(s,t)=W(t)-W(s)$ for $0\leq s\leq t\leq T$ which
produces four more maps from $\{ 0\leq s\leq t\leq T\}$ to $\bbR^D$ where $T>0$ is fixed. We will assume that $\bbS$ and $\bbW$ are connected with $S$
and $W$, respectively, by the, so called, Chen relations which say that for any $0\leq s\leq u\leq t\leq T$,
\begin{eqnarray}\label{2.1}
&\bbS(s,t)=\bbS(s,u)+\bbS(u,t)+S(s,u)\otimes S(u,t)\\
&\mbox{and}\,\,\bbW(s,t)=\bbW(s,u)+\bbW(u,t)+W(s,u)\otimes W(u,t).\nonumber
\end{eqnarray}

For each $\be>0$ and a map $V$ from $\{ 0\leq u\leq v\leq T\}$ to $\bbR^n$ or to $\bbR^n\otimes\bbR^n$ for some $n\geq 1$ we define the $\be$-H\" older
(semi)norm of $V$ on $[s,t],\,0\leq s\leq t\leq T$ by
\[
\| V\|_{\be,[s,t]}=\sup_{s\leq u<v\leq t}\frac {|V(u,v)|}{|v-u|^\be}
\]
where $|\cdot|$ denotes a vector norm in the appropriate Euclidean space. If $[s,t]=[0,T]$ then we will drop the index $[0,T]$ in this notation. Set $R^X(s,t)=X(s,t)-\sig(X(s))S(s,t)$ and $R^\Xi(s,t)=\Xi(s,t)-\sig(\Xi(s))W(s,t)$.
We will assume that for some $\al\in(\frac 13,\frac 12)$,
\begin{equation}\label{2.2}
\| X\|_\al,\,\|\Xi\|_\al,\,\|S\|_\al,\,\| W\|_\al,\,\|\bbS\|_{2\al},\,\|\bbW\|_{2\al},\, \| R^X\|_{2\al}<\infty,\,
\| R^\Xi\|_{2\al}<\infty
\end{equation}
and, also,
\begin{equation}\label{2.3}
\|\sig\|_{C^3},\,\| b\|_{C^2}<\infty
\end{equation}
where $\|\cdot\|_{C^k}$ is the $C^k$-norm.

Next, we introduce
\[
\Psi^X(s,t)=\sig(X(s))S(s,t)+\nabla\sig(X(s))\sig(X(s))\bbS(s,t)+b(X(s))(t-s)
\]
\[
\mbox{and}\,\,\,\Psi^\Xi(s,t)=\sig(\Xi(s))W(s,t)+\nabla\sig(\Xi(s))\sig(\Xi(s))\bbW(s,t)+b(\Xi(s))(t-s)
\]
where
\[
(\nabla\sig(x)\sig(x))_{ijk}=\sum_{l=1}^d\frac {\partial\sig_{ij}(x)}{\partial x_l}\sig_{lk}(x),
\]
\[
(\nabla\sig(x)\sig(x)\bbS(s,t))_i=\sum_{j,k=1}^D(\nabla\sig(x)\sig(x))_{ijk}\bbS^{k,j}(s,t)
\]
and, similarly,
\[
(\nabla\sig(x)\sig(x)\bbW(s,t))_i=\sum_{j,k=1}^D(\nabla\sig(x)\sig(x))_{ijk}\bbW^{k,j}(s,t)
\]
where $\bbS^{k,j}$ and $\bbW^{k,j}$ are $k,j$ components of these matrices.

For any interval $[s,t]\subset[0,T],\, s<t$ define the sequence of partitions $\cP_{s,t}^{(n)}=\{[t^{(n)}_{i-1},t_i^{(n)}],\,
i=1,...,m_n,\, s=t_0^{(n)}<t_1^{(n)}<...<t_{m_n}^{(n)}=t\}$, $m_n=2^n$ of the interval $[s,t]$ inductively by taking
$\cP_{s,t}^{(0)}=\{[s,t]\}$ and building $\cP_{s,t}^{(n+1)}$ from $\cP_{s,t}^{(n)}$ by adding to $t_i^{(n)},\, i=0,1,...,m_n$
the midpoints of all intervals from $\cP_{s,t}^{(n)}$.
\begin{assumption}\label{ass2.1}
Let $\cP_{s,t}^{(n)}=\{[t^{(n)}_{i-1},t_i^{(n)}],\, i=1,...,m_n$ be the sequence of partitions constructed above. Then for any inteval $[s,t]\subset[0,T],\, s<t$,
\begin{equation}\label{2.4}
|X(s,t)-\sum_{i=1}^{m_n}\Psi^X(t^{(n)}_{i-1},t^{(n)}_i)|\to 0\,\,\mbox{as}\,\, n\to\infty
\end{equation}
and
\begin{equation}\label{2.5}
|\Xi(s,t)-\sum_{i=1}^{m_n}\Psi^\Xi(t^{(n)}_{i-1},t^{(n)}_i)|\to 0\,\,\mbox{as}\,\, n\to\infty.
\end{equation}
\end{assumption}
Now we have the following assertion.

\begin{theorem}\label{thm2.2}
Suppose that (\ref{2.1})--(\ref{2.5}) hold true. Then there exists a constant $C_T>0$ which depends only on $T,\,\|\sig\|_{C^3}$ and $\| b\|_{C^2}$ but does not depend on $X,\,\Xi,\, S,\, W,\,\bbS,\,\bbW$ and such that
\begin{equation}\label{2.6}
\| X\|_\al\leq C_T(1+\| S\|_\al+\|\bbS\|_{2\al}),\,\,\| \Xi\|_\al\leq C_T(1+\| W\|_\al+\|\bbW\|_{2\al})
\end{equation}
and
\begin{eqnarray}\label{2.7}
&\| X-\Xi\|_\al\\
&\leq C_T(1+\| S\|_\al+\|\bbS\|_{2\al}+\| W\|_\al+\|\bbW\|_{2\al})^{C_T(1+\| S\|_\al+\|\bbS\|_{2\al}+\| W\|_\al+\|\bbW\|_{2\al})^{1/\al}}\nonumber\\
&\times (\| S-W\|_\al+\sqrt{\|\bbS-\bbW\|_{2\al}}+|X(0)-\Xi(0)|).\nonumber
\end{eqnarray}
\end{theorem}

\subsection{General deterministic setup: discrete time}\label{subsec2.2}

We will deal here with the same maps $\Xi$, $W$, $\bbW$, $\sig$, $b$, as above, but in place of $X$, $S$ and $\bbS$ we will consider now maps
$X_N:[0,T]\to\bbR^d$, $S_N:[0,T]\to\bbR^D$ and $\bbS_N:[0,T]\to\bbR^D\otimes\bbR^D$ such that $X_N(t)=X_N([tN]/N)$, $S_N(t)=S_N([tN]/N)$
and $\bbS_N(s,t)=\bbS_N([sN]/N,[tN]/N)$ assuming that $\bbS_N(k/N,k/N)=0$ for any $k=0,1,...,[tN]$. Hence, $X_N$, $S_N$ and $\bbS_N$ have
now jumps at the points $k/N,\, k=1,2,...,[TN]-1$, and so the standard H\" older semi norms are not appropriate here. By this reason we employ
 here the modified H\" older semi norms defined for each $\be>0$ and a map from $\{ 0\leq u\leq v\leq T\}$ to $\bbR^n$ or to $\bbR^n\otimes\bbR^n$
 (for some $n$) on each interval $[s,t]\subset[0,T],\, s<t$ by
 \[
 \| V\|_{\be,N,[s,t]}=\sup_{s\leq u<v\leq t}\frac {|V(u,v)|}{|u-v|^\be\vee N^{-\be}}
 \]
 where $a\vee b=\max(a,b)$. Again, if $[s,t]=[0,T]$ then we drop the index $[s,t]=[0,T]$ here. We assume that for all $N\geq 1$,
 \begin{equation}\label{2.8}
 \| X_N\|_{\al,N},\|S_N\|_{\al,N},\|\bbS\|_{2\al,N},\,\| R_N\|_{2\al,N}<\infty
 \end{equation}
 where $R_N(s,t)=R^{X_N}(s,t)=X_N(s,t)-\sig(X_N(s))S_N(s,t)$.

 Next, we introduce
 \begin{eqnarray*}
 &\Psi_N^X(s,t)=\sig(X_N(s))S_N(s,t)+\nabla\sig(X_N(s))\sig(X_N(s))\bbS_N(s,t)\\
 &+b(X_N(s))N^{-1}([tN]-[sN])
 \end{eqnarray*}
 where $S_N(s,t)=S_N(t)-S_N(s)$ and the Chen relations
 \begin{equation}\label{2.9}
 \bbS_N(s,t)=\bbS_N(s,u)+\bbS_N(u,t)+S_N(s,u)\otimes S_N(u,t)
 \end{equation}
 are supposed to hold true. We assume again that for each $N\geq 1$,
 \begin{equation}\label{2.10}
 |X_N(s,t)-\sum_{i=1}^{m_n}\Psi_N^X(t_{i-1}^{(n)},t_i^{(n)})|\to 0\,\,\,\mbox{as}\,\,\, n\to\infty
 \end{equation}
 for any sequence of partitions $\cP_{s,t}^{(n)}=\{[t_{i-1}^{(n)},t_i^{(n)}],\, i=1,...,m_n\}$ constructed above. Observe, though, that
 (\ref{2.10}) implies that
 \begin{equation}\label{2.11}
 X_N(s,t)=\sum_{[u,v]\in\cP_{s,t}^{(n)}}\Psi_N^X(u,v)
 \end{equation}
 for all $n\geq n_N=n_N(s,t)=\min\{ k:\, 2^{-k}(t-s)<N^{-1}\}$. Indeed, if $t_i^{(n)}-t_i^{(n-1)}<N^{-1}$ then
 $\Psi_N^X(t_{i-1}^{(n)},t_i^{(n)})=0$ if $k/N\leq t^{(n)}_{i-1}\leq t_i^{(n)}<\frac {k+1}N$ and if
 $k/N\leq t^{(n)}_{i-1}<\frac {k+1}N\leq t_i^{(n)}<\frac {k+2}N$ then
 $\Psi_N^X(t_{i-1}^{(n)},t_i^{(n)})=\Psi_N^X(\frac kN,\frac {k+1}N)$, and so the sum in (\ref{2.10}) remains the same for all
 $n\geq n_N(s,t)$.

 We will assume that for some $\al\in(1/3,1/2)$ (the same as in (\ref{2.8})),
 \begin{equation}\label{2.12}
 \lim_{N\to\infty}N^{-\al}(\| S_N\|_{\al,N}+\sqrt{\|\bbS_N\|_{2\al,N}})=0.
 \end{equation}
 Set also
 \[
 K(\ve)=\min_{K\geq 1}\{ K:\, N^{-\al}(\| S_N\|_{\al,N}+\sqrt{\|\bbS_N\|_{2\al,N}})\leq\ve\,\,\mbox{for all}\,\, N\geq K\}
 \]
 which is finite by (\ref{2.12}) for any $\ve>0$. In this setup we have the following result.

 \begin{theorem}\label{thm2.3} Suppose that $\Xi,\, W,\,\bbW,\sig$ and $b$ satisfy conditions of Theorem \ref{thm2.2}
 while (\ref{2.8})--(\ref{2.10}) and (\ref{2.12}) hold true for $X_N,\, S_N$ and $\bbS_N$. Then there exists a constant
 $C_T>0$ which depends only on $T,\,\|\sig\|_{C^3},\,\| b\|_{C^2}$ and $K(\ve),\,\ve>0$ but does
 not depend on $N,\, X_N,\,\Xi,\, S_N,\, W,\,\bbS_N,\,\bbW$ themselves and such that for all $N\geq 1$,
\begin{equation}\label{2.13}
\| X_N\|_{\al,N}\leq C_T(1+\| S_N\|_{\al,N}+\sqrt{\|\bbS_N\|_{2\al,N}})^{1/\al},\,\,\| \Xi\|_\al\leq C_T(1+\| W\|_\al+\sqrt{\|\bbW\|_{2\al}})^{1/\al}
\end{equation}
and
\begin{eqnarray}\label{2.14}
&\| X_N-\Xi\|_{\al,N}\leq C_T(1+\| S_N\|_{\al,N}+\sqrt{\|\bbS_N\|_{2\al,N}}\\
&+\| W\|_{\al}+\sqrt{\|\bbW\|_{2\al}})^{C_T(1+\| S_N\|_{\al,N}+\sqrt{\|\bbS_N\|_{2\al,N}}+\| W\|_{\al}+\sqrt{\|\bbW\|_{2\al}})^{1/\al}}\nonumber\\
&\times (\| S_N-W\|_{\al,N}+\sqrt{\|\bbS_N-\bbW\|_{2\al,N}}+|X_N(0)-\Xi(0)|)\nonumber
\end{eqnarray}
where $\al\in(1/3,1/2)$ is the same as in (\ref{2.8}) and (\ref{2.12}).
\end{theorem}

This result can be formulated only for the modified H\" older semi norm $\|\cdot\|_{\al,N}$ taking into account that $X_N,\, S_N$ and $\bbS_N$ are
not continuous in time while we still benefit from the fact that their discontinuities are jumps at the specified points of the form $k/N,\, k=1,2,...,
[TN]-1$. By this reason, in order to complement these estimates, we will estimate $X_N-\Xi_N$ below also in the variational norm. For each $p>0$ and a map
$V$ from $\{ 0\leq u\leq v\leq T\}$ to $\bbR^n$ or to $\bbR^n\otimes\bbR^n$ we define the $p$-variational (semi) norm of $V$ on $[0,T]$ by
\[
\| V\|_p=\sup_\cP(\sum_{i=1}^m|V(t_{i-1},t_i)|^p)^{1/p}
\]
where the supremum is taken over all partitions $\cP=\{[t_{i-1},t_i],\, i=1,...,m,\, 0=t_0<t_1<...<t_m=T\}$ of the interval $[0,T]$. To avoid a
confusion between the H\" older and the variational semi norms we will always use the indexes $\al$ or $\be$ for the former and $p$ or $q$ for the
latter.
\begin{corollary}\label{cor2.4}
Suppose that the conditions of Theorem \ref{thm2.3} hold true and $p\al>1,\, p<3$. Then the $p$-variational norm of $X_N-\Xi$ satisfies
\begin{equation}\label{2.15}
\| X_N-\Xi\|_p\leq T^\al\|X_N-\Xi\|_{\al,N}+T^{1/p}N^{-(\al-p^{-1})}(\|\Xi\|_\al+\| X_N\|_{\al,N}).
\end{equation}
\end{corollary}

\subsection{Probabilistic setup: discrete time}\label{subsec2.3}

We will reverse the order and start with the discrete time setup which consists of a complete probability space
$(\Om,\cF,P)$ and a stationary sequence of $d$-dimensional centered random vectors $\xi(n)=(\xi_1(n),...,\xi_d(n))$,
 $-\infty<n<\infty$. We are going to apply the results of Section \ref{subsec2.2} pathwise to
 \begin{eqnarray}\label{2.16}
 &X_N(t)=N^{-1/2}\sum_{0\leq k<[tN]}(\sig(X_N(k/N))\xi(k)+N^{-1}b(X_N(k/N)))\\
 &=\sum_{0\leq k<[tN]}(\sig(X_N(k/N))S_N(\frac kN,\frac {k+1}N)+N^{-1}b(X_N(k/N)))\nonumber
 \end{eqnarray}
where
\[
S_N(t)=N^{-1/2}\sum_{0\leq k<[tN]}\xi(k)\,\,\mbox{and}\,\, S_N(s,t)=S_N(t)-S_N(s),\, s\leq t.
\]
We define also the iterated sums
\[
\bbS_N(s,t)=N^{-1}\sum_{[sN]\leq k<l<[tN]}\xi(k)\otimes\xi(l),\,\,\,\bbS_N(t)=\bbS_N(0,t).
\]

We will rely on strong limit theorems for $S_N$ and $\bbS_N$ from \cite{Ki23+} which require certain moment and weak dependence assumptions
on the process $\xi$. Our setup includes a two parameter family of countably generated $\sig$-algebras
$\cF_{m,n}\subset\cF,\,-\infty\leq m\leq n\leq\infty$ such that
$\cF_{mn}\subset\cF_{m'n'}\subset\cF$ if $m'\leq m\leq n
\leq n'$ where $\cF_{m\infty}=\cup_{n:\, n\geq m}\cF_{mn}$ and $\cF_{-\infty n}=\cup_{m:\, m\leq n}\cF_{mn}$.
It is often convenient to measure the dependence between two sub
$\sig$-algebras $\cG,\cH\subset\cF$ via the quantities
\begin{equation*}
\varpi_{b,a}(\cG,\cH)=\sup\{\| E(g|\cG)-Eg\|_a:\, g\,\,\mbox{is}\,\,
\cH-\mbox{measurable and}\,\,\| g\|_b\leq 1\},
\end{equation*}
where the supremum is taken over real functions and $\|\cdot\|_c$ is the
$L^c(\Om,\cF,P)$-norm. Then more familiar $\al,\rho,\phi$ and $\psi$-mixing
(dependence) coefficients can be expressed via the formulas (see \cite{Bra},
Ch. 4 ),
\begin{eqnarray*}
&\al(\cG,\cH)=\frac 14\varpi_{\infty,1}(\cG,\cH),\,\,\rho(\cG,\cH)=\varpi_{2,2}
(\cG,\cH)\\
&\phi(\cG,\cH)=\frac 12\varpi_{\infty,\infty}(\cG,\cH)\,\,\mbox{and}\,\,
\psi(\cG,\cH)=\varpi_{1,\infty}(\cG,\cH).
\end{eqnarray*}

We set also
\begin{equation*}
\varpi_{b,a}(n)=\sup_{k\geq 0}\varpi_{b,a}(\cF_{-\infty,k},\cF_{k+n,\infty})
\end{equation*}
and accordingly
\[
\al(n)=\frac{1}{4}\varpi_{\infty,1}(n),\,\rho(n)=\varpi_{2,2}(n),\,
\phi(n)=\frac 12\varpi_{\infty,\infty}(n),\, \psi(n)=\varpi_{1,\infty}(n).
\]
We will need also the  approximation rate
\begin{equation*}
\beta (a,l)=\sup_{k\geq 0}\|\xi(k)-E(\xi(k)|\cF_{k-l,k+l})\|_a,\,\, a\geq 1.
\end{equation*}
We will assume that for some $ 1\leq L\leq\infty$, $M$ large enough and $K=\max(2L,4M)$,
\begin{equation}\label{2.17}
E|\xi(0)|^{K}<\infty,\,\,\,\mbox{and}\,\,\,\sum_{k=0}^\infty\sum_{l=k+1}^\infty(\sqrt{\sup_{m\geq l}\be(K,m)}+\varpi_{L,4M}(l))<\infty.
\end{equation}

Next, we introduce the matrix $\Gam=(\Gam_{ij})$ by
\[
\Gam_{ij}=\lim_{k\to\infty}\frac 1k\sum_{n=1}^k\sum_{m=0}^{n-1}E(\xi_i(m)\xi_j(n))=\sum_{l=1}^\infty E(\xi_i(0)\xi_j(l))
\]
and the covariance matrix $\vs=(\vs_{ij})$ by
\begin{equation*}
\vs_{ij}=\lim_{k\to\infty}\frac 1k\sum_{m=0}^k\sum_{n=0}^kE(\xi_i(m)\xi_j(n))=\xi_i(0)\xi_j(0)+\Gam_{ij}+\Gam_{ji}
\end{equation*}
where the limits are known to exist under our conditions (see \cite{Ki23}). Let $\cW=(\cW^1,\cW^2,...,\cW^d)$ be a $d$-dimensional
Brownian motion with the covariance matrix $\vs$ (at the time 1) and introduce the rescaled Brownian motion $W_N(t)=N^{-1/2}\cW(Nt),\,
t\in[0,T],\, N\geq 1$. Let $\Xi=\Xi_N$ be the diffusion solving the stochastic differential equation
\begin{equation}\label{2.18}
d\Xi_N(t)=\sig(\Xi_N(t))dW_N(t)+(b(\Xi_N(t))+c(\Xi_N(t))dt
\end{equation}
where
\[
c(x)=\nabla\sig(x)\sig(x)\Gam\,\,\mbox{and in coordinates}\,\, c_i(x)=\sum_{j,k}(\nabla\sig(x)\sig(x))_{ijk}\Gam_{kj}.
\]

Next, we introduce
\begin{equation}\label{2.19}
\bbW_N(s,t)=\int_s^tW_N(s,v)\otimes dW_N(v)+(t-s)\Gam
\end{equation}
which can be written in the coordinate-wise form as
\[
\bbW_N^{ij}(s,t)=\int_s^tW_N^i(s,v)dW^j_N(v)+(t-s)\Gam_{ij}.
\]
We will show relying on Theorem \ref{thm2.3}, Corollary \ref{cor2.4} and the results from \cite{Ki23} and \cite{Ki23+}
that the following assertion holds true.

\begin{theorem}\label{thm2.5} Suppose that (\ref{2.17}) is satisfied. Then the stationary sequence of random vectors $\xi(n),\,-\infty<n<\infty$
 can be redefined preserving its distributions on a sufficiently rich probability space which contains also a $d$-dimensional Brownian motion
 $\cW$ with the covariance matrix $\vs$ (at the time 1) so that the processes $X_N$ and $\Xi_N$ constructed by $\xi$ and the rescaled processes
 $W_N(t)=N^{-1/2}\cW(Nt)$ satisfy
 \begin{equation}\label{2.20}
 \| X_N-\Xi_N\|_p=O(N^{-\del})\quad\mbox{a.s.}
 \end{equation}
 provided $X_N(0)=\Xi_N(0)$, where $p\in(2,3)$ and $\del>0$ does not depend on $N$.
 \end{theorem}

  Important classes of processes satisfying our conditions come from
dynamical systems. Let $F$ be a $C^2$ Axiom A diffeomorphism (in
particular, Anosov) in a neighborhood $\Om$ of an attractor or let $F$ be
an expanding $C^2$ endomorphism of a compact Riemannian manifold $\Om$ (see
\cite{Bow}), $g$ be either a H\" older continuous vector function or a
vector function which is constant on elements of a Markov partition and let $\xi(n)=
\xi(n,\om)=g(F^n\om)$. Here the probability space is $(\Om,\cB,P)$ where $P$ is a Gibbs
 invariant measure corresponding to some H\"older continuous function and $\cB$ is the Borel $\sig$-field.
  Let $\zeta$ be a finite Markov partition for $F$ then we can take $\cF_{kl}$
 to be the finite $\sig$-algebra generated by the partition $\cap_{i=k}^lF^i\zeta$.
 In fact, we can take here not only H\" older continuous $g$'s but also indicators
of sets from $\cF_{kl}$. Our conditions allow all such functions
since the dependence of H\" older continuous functions on $m$-tails, i.e. on events measurable
with respect to $\cF_{-\infty,-m}$ or $\cF_{m,\infty}$, decays exponentially fast in $m$ while
our conditions are much weaker than that. A related class of dynamical systems
corresponds to $F$ being a topologically mixing subshift of finite type which means that $F$
is the left shift on a subspace $\Om$ of the space of one (or two) sided
sequences $\om=(\om_i,\, i\geq 0), \om_i=1,...,l_0$ such that $\om\in\Om$
if $\pi_{\om_i\om_{i+1}}=1$ for all $i\geq 0$ where $\Pi=(\pi_{ij})$
is an $l_0\times l_0$ matrix with $0$ and $1$ entries and such that $\Pi^n$
for some $n$ is a matrix with positive entries. Again, we have to take in this
case $g$ to be a H\" older continuous bounded function on the sequence space above,
 $P$ to be a Gibbs invariant measure corresponding to some H\" older continuous function and to define
$\cF_{kl}$ as the finite $\sig$-algebra generated by cylinder sets
with fixed coordinates having numbers from $k$ to $l$. The
exponentially fast $\psi$-mixing is well known in the above cases (see \cite{Bow}) and this property
is much stronger than what we assume. Among other
dynamical systems with exponentially fast $\psi$-mixing we can mention also the Gauss map
$Fx=\{1/x\}$ (where $\{\cdot\}$ denotes the fractional part) of the
unit interval with respect to the Gauss measure $G$ and more general transformations generated
by $f$-expansions (see \cite{Hei}). Gibbs-Markov maps which are known to be exponentially fast
$\phi$-mixing (see, for instance, \cite{MN}) can be also taken as $F$ in $\xi(n)=g\circ F^n$ as above.

\subsection{Probabilistic continuous time setup}\label{subsec2.4}

Here we obtain $X_N$ as the solution of the equation
\begin{eqnarray}\label{2.21}
&X_N(t)=X_N(s)+N^{1/2}\int_s^t\sig(X_N(u))\xi(uN)du+\int_s^tb(X_N(u))du\\
&=X_N(s)+\int_s^t\sig(X_N(u))dS_N(u)+\int_s^tb(X_N(u))du\nonumber
\end{eqnarray}
where $S_N(u)=N^{-1/2}\int_0^{uN}\xi(v)dv$, $\xi$ is a centered vector stationary stochastic process in $\bbR^d$ and we assume
that this integral exists. As above, we suppose
 that $\sig$ and $b$ have bounded $C^3$ and $C^2$ norms, respectively, and under the moment condition (\ref{2.17}) on $\xi$
it is not difficult to see, say, by the Picard approximations or by contracting mappings, that for $P$-almost all $\om\in\Om$ there exists
a unique solution $X_N(t)=X_N(t,\om),\, t\in[0,T]$ of the integral equation (\ref{2.21}). Next, we set
\[
\bbS_N(s,t)=\int_s^tS_N(s,u)dS_N(u)=\int_s^tS_N(s,u)\xi(u)du
\]
assuming this integral exists where, as before, $S_N(s,u)=S_N(u)-S_N(s)$. In order to employ Theorem \ref{thm2.2} we will assume also
that a.s. for all $N$,
\begin{equation}\label{2.22}
\| S_N\|_\al<\infty\,\,\,\mbox{and}\,\,\,\|\bbS_N\|_{2\al}<\infty
\end{equation}
while $\| X_N\|_\al<\infty$ and $\| R_N\|_{2\al}<\infty$ we obtain here automatically assuming, for instance,
that $\xi$ is a continuous in time process. 

In the continuous time case we will consider two types of assumptions concerning the process $\xi$.

\subsubsection{Straightforward setup}\label{subsubsec2.4.1}
Our direct continuous time setup consists of a centered $d$-dimensional stationary process
$\xi(t),\, t\geq 0$ on a probability space $(\Om,\cF,P)$ and of a family of
$\sig$-algebras $\cF_{st}\subset\cF,\,-\infty\leq s\leq t\leq\infty$ such
that $\cF_{st}\subset\cF_{s't'}$ if $s'\leq s$ and $t'\geq t$. For all $t\geq 0$ we set
\begin{equation*}
\varpi_{b,a}(t)=\sup_{s\geq 0}\varpi_{b,a}(\cF_{-\infty,s},\cF_{s+t,\infty})
\end{equation*}
and
\begin{equation*}
\beta (a,t)=\sup_{s\geq 0}\|\xi(s)-E(\xi(s)|\cF_{s-t,s+t})\|_a.
\end{equation*}
where $\varpi_{b,a}(\cG,\cH)$ is defined in Section \ref{subsec2.3}.  We continue to
impose the assumption (\ref{2.17}) on the  decay rates  of
$\varpi_{b,a}(t)$ and $\beta (a,t)$. Although they only involve integer
values of $t$, it will  suffice since these are non-increasing functions of $t$.

Next, we introduce the covariance matrix $\vs=(\vs_{ij})$ defined by
\begin{equation*}
\vs_{ij}=\lim_{t\to\infty}\frac 1t\int_{0}^t\int_{0}^t\vs_{ij}(u-v)dudv,\,\,
\mbox{where}\,\, \vs_{ij}(u-v)=E(\xi_i(u)\xi_j(v))
\end{equation*}
and the limit here exists under our conditions in the same way as in the discrete time setup.
In order to formulate our results we define
\begin{eqnarray*}
&S_N(t)=N^{-1/2}\int_0^t\xi(u)du,\,\, S_N(s,t)=S_N(t)-S_N(s)\\
&\mbox{and}\,\,\,\bbS_N(s,t)=N^{-1}\int_s^t\xi(u)du\int_s^u\xi(v)dv,\,\,\bbS_N(t)=\bbS_N(0,t).
\end{eqnarray*}
We will deal here with the Brownian motion $W_N$ with the covariance matrix $\vs$ defined above and with
 $\bbW_N$ defined by (\ref{2.19}) with the matrix $\Gam$ given by
\[
\Gam=(\Gam_{ij}),\,\,\Gam_{ij}=\int_0^\infty E(\xi_i(0)\xi_j(u))du+\int_0^1du\int_0^uE(\xi_i(v)\xi_j(u))dv
\]

\subsubsection{Suspension setup}\label{subsubsec2.4.2}
Here we start with a complete probability space $(\Om,\cF,P)$, a
$P$-preserving invertible transformation $\vt:\,\Om\to\Om$ and
a two parameter family of countably generated $\sig$-algebras
$\cF_{m,n}\subset\cF,\,-\infty\leq m\leq n\leq\infty$ such that
$\cF_{mn}\subset\cF_{m'n'}\subset\cF$ if $m'\leq m\leq n
\leq n'$ where $\cF_{m\infty}=\cup_{n:\, n\geq m}\cF_{mn}$ and
$\cF_{-\infty n}=\cup_{m:\, m\leq n}\cF_{mn}$. The setup includes
also a (roof or ceiling) function $\tau:\,\Om\to (0,\infty)$ such that
for some $\hat L>0$,
\begin{equation*}
\hat L^{-1}\leq\tau\leq\hat L.
\end{equation*}
Next, we consider the probability space $(\hat\Om,\hat\cF,\hat P)$ such that $\hat\Om=\{\hat\om=
(\om,t):\,\om\in\Om,\, 0\leq t\leq\tau(\om)\},\, (\om,\tau(\om))=(\vt\om,0)\}$, $\hat\cF$ is the
restriction to $\hat\Om$ of $\cF\times\cB_{[0,\hat L]}$, where $\cB_{[0,\hat L]}$ is the Borel
$\sig$-algebra on $[0,\hat L]$ completed by the Lebesgue zero sets, and for any $\Gam\in\hat\cF$,
\[
\hat P(\Gam)=\bar\tau^{-1}\int\bbI_\Gam(\om,t)dtdP(\om)\,\,\mbox{where}\,\,\bar\tau=\int\tau dP=E\tau
\]
and $E$ denotes the expectation on the space $(\Om,\cF,P)$.

Finally, we introduce a vector valued stochastic process $\xi(t)=\xi(t,(\om,s))$, $-\infty<t<\infty,\, 0\leq
s\leq\tau(\om)$ on $\hat\Om$ satisfying
\begin{eqnarray*}
&\int\xi(t)d\hat P=0,\,\xi(t,(\om,s))=\xi(t+s,(\om,0))=\xi(0,(\om,t+s))\,\,\mbox{if}\,\, 0\leq t+s<\tau(\om)\\
&\mbox{and}\,\,
\xi(t,(\om,s))=\xi(0,(\vt^k\om,u))\,\,\mbox{if}\,\, t+s=u+\sum_{j=0}^{k-1}\tau(\vt^j\om)\,\,\mbox{and}\,\,
0\leq u<\tau(\vt^k\om).
\end{eqnarray*}
This construction is called in dynamical systems a suspension and it is a standard fact that $\xi$ is a
stationary process on the probability space $(\hat\Om,\hat\cF,\hat P)$ and in what follows we will write
also $\xi(t,\om)$ for $\xi(t,(\om,0))$. In this setup we define $S_N$ and $\bbS_N$ by normalized integrals
and iterated integrals, as above.

Set $\eta(\om)=\int_0^{\tau(\om)}\xi(s,\om)ds$, $\eta(m)=\eta\circ\vt^m$ and
\begin{equation*}
\be(a,l)=\sup_m\max\big(\|\tau\circ\vt^m-E(\tau\circ\vt^m|\cF_{m-l,m+l})\|_a,\,\|\eta(m)-E(\eta(m)|\cF_{m-l,m+l})\|_a\big).
\end{equation*}
We define $\varpi_{b,a}(n)$ by (\ref{2.2}) with respect to the $\sig$-algebras $\cF_{kl}$ appearing here.
Observe also that $\eta(k)=\eta\circ\vt^k$ is a stationary sequence of random vectors and we introduce also the covariance matrix
 \begin{equation*}
\vs_{ij}=\lim_{n\to\infty}\frac 1n\sum_{k,l=0}^nE(\eta_i(k)\eta_j(l))
\end{equation*}
where the limit exists under our conditions in the same way as before.
We will deal here with the Brownian motion $W_N$ with the covariance matrix $\vs$ defined above and with
 $\bbW_N$ defined by (\ref{2.19}) with the matrix $\Gam$ given by
\[
\Gam=(\Gam_{ij}),\,\,\Gam_{ij}=\sum_{l=1}^\infty E(\eta_i(0)\eta_j(l))+E\int_0^{\tau(\om)}\xi_j(s,\om)ds\int_0^s\xi_i(u,\om)du.
\]

In the continuous time case we have the following assertion.
\begin{theorem}\label{thm2.6} Suppose that (\ref{2.17}) and (\ref{2.22}) are satisfied. Then in both the straightforward and suspension setups
the stationary process $\xi(n),\,-\infty<n<\infty$ can be redefined preserving its distributions on a sufficiently rich probability
 space which contains also a $d$-dimensional Brownian motion  $\cW$ with the covariance matrix $\vs$ (at the time 1) so that the
 process $\Xi$ solving (\ref{2.18}) with the rescaled processes $W_N(t)=N^{-1/2}\cW(Nt)$ and the matrix $\Gam$ together with $X_N$
  constructed by $\xi$ satisfy
 \begin{equation}\label{2.23}
 \| X_N-\Xi_N\|_p\leq T^{\al}\| X_N-\Xi_N\|_\al=O(N^{-\del})\quad\mbox{a.s.}
 \end{equation}
 provided $X_N(0)=\Xi_N(0)$, where $p\in(2,3)$, $\al=1/p$ and $\del>0$ does not depend on $N$.
 \end{theorem}
 This theorem extends applicability of our results to processes $\xi(t)=g\circ F^t$ constructed by a hyperbolic flow $F^t$ (see \cite{BR})
 on a compact Riemannian manifold $\hat\Om$ preserving a Gibbs measure $\hat P$ built by a H\" older continuous function with $g$ being H\" older,
 as well.

\begin{remark}\label{rem2.7}
In fact, our deterministic setup allows also slightly more general forms of the processes $X_N$ considered above. Namely, in the discrete time setup let
\begin{eqnarray*}
&X_N(t)=N^{-1/2}\sum_{0\leq k<[tN]}(\tilde\sig(X_N(k/N))\tilde\xi(k)+\tilde{\tilde\sig}(X_N(k/N))\tilde{\tilde\xi}(k))\\
&=\sum_{0\leq k<[tN]}\sig(X_N(k/N))S_N(\frac kN,\frac {k+1}N)
\end{eqnarray*}
where $\sig_{ij}(x)=\tilde\sig_{ij}(x)$ and $\sig_{i,j+d}(x)=\tilde{\tilde\sig}_{ij}(x)$ when $i,j=1,...,d$, i.e. $\sig$ is now a $d\times 2d$-matrix,
and $S^i_N(t)=N^{-1/2}\sum_{0\leq k<[tN]}\tilde\xi_i(k)$, $S^{i+d}_N(t)=N^{-1}\sum_{0\leq k<[tN]}\tilde{\tilde\xi}_i(k)$ for $i=1,...,d$. In
the continuous time case the definitions are similar with integrals in place of sums. The same proof as in this paper will yield the almost sure
estimates of the form (\ref{2.20}).
\end{remark}

\begin{remark}\label{rem2.8}
In the continuous time case we can deal also with the more frequent form of (\ref{1.1}) when $N$ is replaced by $\ve^{-2}$ so that we have
\begin{equation}\label{2.24}
\frac {dX^\ve(t)}{dt}=\ve^{-1}\sig(X^\ve(t))\xi(t/\ve^2)+b(X^\ve(t))
\end{equation}
and the goal is to have (\ref{2.20}) of Theorem \ref{thm2.5} in the form
\begin{equation}\label{2.25}
\| X^\ve-\Xi^\ve\|_\al=O(\ve^\del)\quad\mbox{a.s.}\,\,\,\mbox{as}\,\,\,\ve\to 0
\end{equation}
where the diffusion $\Xi^\ve$ solves the stochastic differential equation (\ref{2.18}) with $W_N$ replaced by $W^\ve(t)=\ve\cW(\ve^{-2}t),
\, 0\leq t\leq T$. We proved that (\ref{2.25}) holds true when $\ve=\ve_N=N^{-1/2}$ and it turns out that this implies (\ref{2.25}) as
$\ve\to 0$ in general. Indeed, our deterministic setup Theorem \ref{thm2.2} does not depend on parameters either $N$ or $\ve^{-2}$, and
so we can rely on its estimate (\ref{2.7}). Hence, we only have to check that the asymptotic behavior as $\ve\to 0$ of
\[
S^\ve(t)=\ve\int_0^{t\ve^{-2}}\xi(u)du\quad\mbox{and of}\quad\bbS^\ve(t)=\ve^2\int_0^{t\ve^{-2}}\xi(u)du\int_0^u\xi(v)dv,\,\, t\in[0,T]
\]
is, essentially, the same as of
\begin{eqnarray*}
&S_{[\ve^{-2}]}(t)=[\ve^{-2}]^{-1/2}\int_0^{t[\ve^{-2}]}\xi(u)du\quad\mbox{and of}\\
&\bbS_{[\ve^{-2}]}(t)=[\ve^{-2}]^{-1}\int_0^{t[\ve^{-2}]}\xi(u)du\int_0^u\xi(v)dv,\,\, t\in[0,T],
\end{eqnarray*}
respectively. We have
\begin{eqnarray*}
&\sup_{0\leq t\leq T}|S^\ve(t)-S_{[\ve^{-2}]}(t)|\leq |\ve-[\ve^{-2}]^{-1/2}|\sup_{0\leq t\leq T/\ve^2}|\int_0^t\xi(u)du|\\
&+[\ve^{-2}]^{-1/2}\sup_{0\leq n\leq [\ve^{-2}]}\int_{nT}^{(n+1)T}|\xi(u)|du=O(\ve^\del)\quad\mbox{a.s.}
\end{eqnarray*}
and, similarly,
\[
\sup_{0\leq t\leq T}|\bbS^\ve(t)-\bbS_{[\ve^{-2}]}(t)|=O(\ve^\del)\quad\mbox{a.s.}
\]
for some $\del>0$. Using an argument similar to the one employed in Section \ref{subsec5.1} below we can extend these estimates to a
$\al$-H\" older semi norm with $\al\in(1/3,1/2)$. This together with the estimate (\ref{2.7}) of Theorem \ref{thm2.2} yields (\ref{2.25})
similarly to Theorem \ref{thm2.6}.
\end{remark}

\section{Bounds on $X$ and $\Xi$ in the deterministic setup }\label{sec3}\setcounter{equation}{0}

\subsection{Continuous time case}\label{subsec3.1}
 We will partially follow here Chapter 4 in \cite{FH}.
Let $\cP_{s,t}^{(n)}$ be a sequence of partitions of $[s,t]$
defined in Section \ref{subsec2.1}. Set
\[
\del\Psi^X(s,u,t)=\Psi^X(s,t)-\Psi^X(s,u)-\Psi^X(u,t)
\]
and
\[
I^n(s,t)=\sum_{[u,v]\in\cP_{s,t}^{(n)}}\Psi^X(u,v)
\]
for $n\geq 1$, taking $I^0(s,t)=\Psi^X(s,t)$. Then
\[
I^{n+1}(s,t)=I^n(s,t)-\sum_{[u,v]\in\cP_{s,t}^{(n)}}\del\Psi^X(u,m,v)
\]
where $m$ is the midpoint of the interval $[u,v]$, and so
\[
|I^{n+1}(s,t)-I^n(s,t)|=|\sum_{[u,v]\in\cP_{s,t}^{(n)}}\del\Psi^X(u,m,v)|.
\]
Since by Assumption \ref{ass2.1}, $I^n(s,t)\to X(s,t)$ as $n\to\infty$ we obtain that
\begin{equation}\label{3.1}
|X(s,t)-\Psi^X(s,t)|\leq\sum_{n=0}^\infty|\sum_{[u,v]\in\cP_{s,t}^{(n)}}\del\Psi^X(u,m,v)|.
\end{equation}

Let $\be=3\al$ with $\frac 13<\al<\frac 12$ and set
\[
\|\del\Psi^X\|_{\be,[s,t]}=\sup_{s\leq u<w<v\leq t}\frac {|\del\Psi^X(u,w,v)|}{|v-u|^\be},\,\,\,[s,t]\subset[0,T].
\]
If $[u,v]\in\cP_{s,t}^{(n)}$ then $v-u=2^{-n}(t-s)$, and so for any $m\in[u,v]$,
\[
|\del\Psi^X(u,m,v)|\leq \|\del\Psi^X\|_{\be,[u,v]}|v-u|^\be\leq 2^{-n\be}\|\del\Psi^X\|_{\be,[s,t]}|t-s|^\be.
\]
This together with (\ref{3.1}) and the fact that $\cP_{s,t}^{(n)}$ contains $2^n$ intervals yields the estimate
\begin{equation}\label{3.2}
|X(s,t)-\Psi^X(s,t)|\leq (1-2^{1-\be})^{-1}\|\del\Psi^X\|_{\be,[s,t]}|t-s|^\beta
\end{equation}
taking into account that $\be=3\al>1$.

 Next, we write $\del\Psi^X$ in a more convenient form. Relying on (\ref{2.1}) we have
 \begin{eqnarray}\label{3.3}
 &\,\,\,\del\Psi^X(s,u,t)=\sig(X(s))S(s,t)+\nabla\sig(X(s))\sig(X(s))\bbS(s,t)+b(X(s))(t-s)\\
 &-\sig(X(s))S(s,u)-\nabla\sig(X(s))\sig(X(s))\bbS(s,u)-b(X(s))(u-s)\nonumber\\
 &-\sig(X(u))S(u,t)-\nabla\sig(X(u))\sig(X(u))\bbS(u,t)-b(X(u))(t-u)\nonumber\\
 &=-R^{\sig(X)}(s,u)S(u,t)-\nabla\sig(X)\sig(X)(s,u)\bbS(u,t)-b(X)(s,u)(t-u)\nonumber
 \end{eqnarray}
 where $b(X)(s,u)=b(X(u))-b(X(s))$, $\nabla\sig(X)\sig(X)(s,u)=\nabla\sig(X(u))\sig(X(u))-\nabla\sig(X(s))\sig(X(s))$,
 $R^{\sig(X)}(s,u)=\sig(X)(s,u)-\nabla\sig(X(s))X(s,u)+\nabla\sig(X(s))R^X(s,u)$, $\sig(X)(s,u)=\sig(X(u))-\sig(X(s))$,
 $(\nabla\sig(x)y)_{ij}=\sum_{k=1}^d\frac {\partial\sig_{ij}(x)}{\partial x_k}y_k$ and, recall,
  $R^X(s,u)=X(s,u)-\sig(X(s))S(s,u)$.

 By (\ref{3.2}) with $\be=3\al$,
 \begin{eqnarray}\label{3.4}
&|R^X(s,t)|\leq|X(s,t)-\Psi^X(s,t)|+|\nabla\sig(X(s))\sig(X(s))\bbS(s,t)|\\
&+|b(X(s))|(t-s)\leq(1-2^{1-3\al})^{-1}\|\del\Psi^X\|_{3\al,[s,t]}|t-s|^{3\al}\nonumber\\
&+\|\nabla\sig\|_\infty\|\sig\|_\infty\|\bbS\|_{2\al,[s,t]}|t-s|^{2\al}+\| b\|_\infty|t-s|\nonumber
\end{eqnarray}
where $\|\cdot\|_\infty$ is the supremum (or $L^\infty$) norm. Combining (\ref{3.3}) and (\ref{3.4}) we obtain
\begin{eqnarray}\label{3.5}
&|R^X(s,t)|\leq(1-2^{1-3\al})^{-1}(\| R^{\sig(X)}\|_{2\al,[s,t]}\| S\|_{\al,[s,t]}\\
&+\|\nabla\sig(X)\sig(X)\|_{\al,[s,t]}\|\bbS\|_{2\al,[s,t]}+\| b(X)\|_{\al,[s,t]}(t-s)^{1-2\al})(t-s)^{3\al})\nonumber\\
&+\|\nabla\sig\|_\infty\|\sig\|_\infty\|\bbS\|_{2\al,[s,t]}|t-s|^{2\al}+\| b\|_\infty|t-s|.\nonumber
\end{eqnarray}
Since by the two term Taylor formula with the remainder
\[
|\sig(X)(s,t)-\nabla\sig(X(s))X(s,t)|\leq\frac 12\|\nabla^2\sig\|_\infty|X(s,t)|^2,
\]
we see by the formula for $R^{\sig(X)}$ above that
\begin{equation}\label{3.6}
\| R^{\sig(X)}\|_{2\al,[s,t]}\leq\frac 12\|\nabla^2\sig\|_\infty\| X\|^2_{\al,[s,t]}+\|\nabla\sig\|_\infty\| R^X\|_{2\al,[s,t]}.
\end{equation}
This together with (\ref{2.2}) and (\ref{3.3}) yield that
\begin{equation}\label{3.6+}
\|\del\Psi^X\|_\be<\infty
\end{equation}
which, in particular, makes the estimate (\ref{3.2}) useful.

Next, for each $h>0$ we will use the semi norms
\[
\|\cdot\|_{\al,h}=\sup_{I\subset[0,T],\,|I|\leq h}\|\cdot\|_{\al,I}\,\,\mbox{and}\,\,\|\cdot\|_{2\al,h}=\sup_{I\subset[0,T],\,|I|\leq h}\|\cdot\|_{2\al,I}.
\]
Then, (\ref{3.5}) and (\ref{3.6}) imply
\begin{eqnarray}\label{3.7}
&\|R^X\|_{2\al,h}\leq(1-2^{1-3\al})^{-1}(\| R^{\sig(X)}\|_{2\al,h}\| S\|_{\al,h}\\
&+\|\nabla\sig(X)\sig(X)\|_{\al,h}\|\bbS\|_{2\al,h}+\| b(X)\|_{\al,h}h^{1-2\al})h^\al+\|\nabla\sig\|_\infty\|\sig\|_\infty\|\bbS\|_{2\al,h}\nonumber\\
&+\| b\|_\infty h^{1-2\al}\leq(1-2^{1-3\al})^{-1}h^\al(\frac 12\|\nabla^2\sig\|_\infty\| X\|^2_{\al,h}\| S\|_{\al,h}\nonumber\\
&+\|\nabla\sig\|_\infty\| R^X\|_{2\al,h}\| S\|_{\al,h}+(\|\nabla\sig\|_\infty^2+\|\nabla^2\sig\|_\infty\|\sig\|_\infty)\| X\|_{\al,h}\|\bbS_N\|_{2\al,h}\nonumber\\
&+\|\nabla b\|_\infty\| X\|_{\al,h}h^{1-2\al})+\|\nabla\sig\|_\infty\|\sig\|_\infty\|\bbS\|_{2\al,h}+\| b\|_\infty h^{1-2\al}\nonumber
\end{eqnarray}
where we use that
\begin{eqnarray*}
&\|\nabla\sig(X)\sig(X)\|_{\al,h}=\sup_{[u,v]\subset[0,T],0<v-u\leq h}\frac {|\nabla\sig(X)\sig(X)(u,v)|}{|v-u|^\al}\\
&\leq(\|\nabla\sig\|^2_\infty+\nabla^2\sig\|_\infty\|\sig\|_\infty)\| X\|_{\al,h}.
\end{eqnarray*}

Set
\[
C_\sig=2(\|\nabla\sig\|_\infty+\|\nabla^2\sig\|_\infty+\|\nabla\sig\|^2_\infty+\|\nabla^2\sig\|_\infty\|\sig\|_\infty)
\]
and
\[
h_0=\sup\{ h\leq 1:\, h^\al C_\al(\| S\|_\al+\sqrt{\|\bbS\|_{2\al}})+h^{1-2\al}(\| b\|_\infty+h^\al\|\nabla b\|_\infty)\leq 1-2^{1-3\al}\}.
\]
Then (\ref{3.7}) implies that for any $h\leq h_0$,
\begin{eqnarray*}
&\| R^X\|_{2\al,h}\leq\frac 12\| X\|^2_{\al,h}+\frac 12\| R^X\|_{2\al,h}+\| X\|_{\al,h}\sqrt{\|\bbS\|_{2\al,h}}\\
&+\| X\|_{\al,h}+\|\nabla\sig\|_\infty\|\sig\|_\infty\|\bbS\|_{2\al,h}+1,
\end{eqnarray*}
and so
\begin{eqnarray}\label{3.8}
&\| R^X\|_{2\al,h}\leq\| X\|^2_{\al,h}+2\| X\|_{\al, h}\sqrt{\|\bbS\|_{2\al,h}}+2\|\nabla\sig\|_\infty\|\sig\|_\infty\|\bbS\|_{2\al,h}\\
&+2\| X\|_{\al,h}+2\leq 2\| X\|^2_{\al,h}+(2\|\nabla\sig\|_\infty\|\sig\|_\infty+1)\|\bbS\|_{2\al,h}+2\| X\|_{\al,h}+2.\nonumber
\end{eqnarray}

On the other hand,
\[
X(s,t)=\sig(X(s))S(s,t)+R^X(s,t),\,\,\mbox{and so}\,\,\| X\|_{\al,h}\leq\|\sig\|_\infty\| S\|_{\al,h}+\| R^X\|_{2\al,h}h^\al.
\]
This together with (\ref{3.8}) yields that
\begin{eqnarray}\label{3.9}
&\| X\|_{\al,h}\leq\|\sig\|_\infty\| S\|_{\al,h}+2h^\al\| X\|^2_{\al,h}\\
&+h^\al(2\|\nabla\sig\|_\infty\|\sig\|_\infty+1)\|\bbS\|_{2\al,h}+2h^\al\| X\|_{\al,h}+2h^\al.\nonumber
\end{eqnarray}
Multiplying both parts of (\ref{3.9}) by $4h^\al$ and setting $\vf_h=4h^\al\| X\|_{\al,h}$ and
\[
\la_h=4h^\al\|\sig\|_\infty\| S\|_\al+4h^{2\al}(2\|\nabla\sig\|_\infty\|\sig\|_\infty+1)\|\bbS\|_{2\al,h}+8h^{2\al}
\]
we obtain that
\begin{equation}\label{3.10}
\vf_h\leq\la_h+\frac 12\vf_h^2+2h^\al\vf_h.
\end{equation}

Set
\[
h_1=\max\{ h\leq h_0:\,\la_h\leq\frac 5{72}\}.
\]
Since $8h^{2\al}\leq\frac 5{72}$ then $2h^\al<1/4$ and (\ref{3.10}) implies
\begin{equation}\label{3.11}
\vf_h\leq 2\la_h+\vf_h^2.
\end{equation}
Since $\la_h$ is increasing in $h$ we will have $2\la_h\leq\frac 5{36}$ for all $h\leq h_1$, and so for such $h$ either
\[
\vf_h\geq\vf_+=\frac 12+\sqrt{\frac 14-2\la_h}\geq\frac 56\,\,\,\mbox{or}
\]
\[
\vf_h\leq\vf_{-}=\frac 12-\sqrt{\frac 14-2\la_h}=\frac 12(1-\sqrt{1-8\la_h})\leq\frac 16.
\]

Since $\| X\|_\al<\infty$ by the assumption (\ref{2.2}), $\vf_h\to 0$ as $h\to 0$. Observe also that
\begin{eqnarray*}
&\| X\|_{\al,h}=\sup_{[u,v]\subset[0,T],\, 0<v-u\leq h}\frac {|X(v)-X(u)|}{|v-u|^\al}\\
&\leq\sup_{[u,v]\subset[0,T],\, 0<v-u\leq h}\frac {|X(v)-X(v-\frac {v-u}2)|+|X(v-\frac {v-u}2)-X(u)|}{|v-u|^\al}\leq 2\| X\|_{\al,\frac h2}.
\end{eqnarray*}
Hence, $\vf_h\leq 2^{1+\al}\vf_{h/2}<3\vf_{h/2}$, and so by the monotonicity of $\vf_h$ in $h$ it cannot jump from below $\vf_{-}$ to
above $\vf_+$ which means that $\vf_h\leq\vf_{-}$ for all $h\leq h_1$. Thus, $\vf_h\leq 1/6$ when $h\leq h_1$ and then $\frac 56\vf_h\leq 2\la_h$,
i.e. $\vf_h\leq\frac {12}5\la_h$. It follows that for all $h\leq h_1$,
\[
\| X\|_{\al,h}\leq\frac {12}5(\|\sig\|_\infty\| S\|_\al+h^\al(2\|\nabla\sig\|_\infty\|\sig\|_\infty+1)\|\bbS\|_{2\al}+2h^\al).
\]
Since $4h^{2\al}\|\bbS\|_{2\al}<\frac 5{72}$ and $8h^{2\al}<\frac 5{72}$ when $h\leq h_1$, we see that $\frac {12}5\sqrt{\|\bbS\|_{2\al}}<1$
and $\frac {24}5h^\al<1$, and so the above estimate can be rewritten as
\begin{equation}\label{3.12}
\| X\|_{\al,h}\leq\tilde C_{\sig}(\| S\|_\al+\sqrt{\|\bbS\|_{2\al}}+1)
\end{equation}
where
\[
\tilde C_\sig=\frac {12}5\|\sig\|_\infty+2\|\nabla\sig\|_\infty\|\sig\|_\infty+1.
\]

Finally,
\begin{eqnarray*}
&\| X\|_\al=\sup_{0\leq s<t\leq T}\frac {|X(s,t)}{|t-s|^\al}\\
&\leq\sup_{0\leq s<t\leq T}\frac {\sum_{k=0}^{\nu(s,t)-1}|X(s+kh_1,s+(k+1)h_1)|+|X(s+\nu(s,t)h_1,t)|}{(\nu(s,t))^\al h_1^\al+(t-s-\nu(s,t)h_1)^\al}\\
&\leq\| X\|_{\al,h_1}(1+\sup_{0\leq s<t\leq T}(\nu(s,t))^{1-\al})=\| X\|_{\al,h_1}(1+[T/h_1]^{1-\al})\\
\end{eqnarray*}
where $\nu(s,t)=[\frac {t-s}{h_1}]$. This together with (\ref{3.12}) and the definition of $h_0,\, h_1$ and $\la_h$ yields the first half of (\ref{2.6})
in Theorem \ref{thm2.2}. The second half of (\ref{2.6}) follows by exactly the same proof with $X,\, S$ and $\bbS$ replaced by $\Xi,\, W$ and $\bbW$
as there are no differences between our assumptions concerning these two sets of maps.

\subsection{Discrete time case}\label{subsec3.2}
We set
\[
\del\Psi^X_N(s,u,t)=\Psi_N^X(s,t)-\Psi^X_N(s,u)-\Psi_N^X(u,t)
\]
and
\[
\|\del\Psi_N^X\|_{\be,N,[s,t]}=\sup_{s\leq u<w<v\leq t}\frac {|\del\Psi^X_N(u,w,v)|}{|v-u|^\be\vee N^{-\be}}
\]
for $[s,t]\subset[0,T],\,\be=3\al$ and $1/3<\al<1/2$. Then
\[
|\del\Psi_N^X(u,m,v)|\leq\|\del\Psi_N^X\|_{\be,N,[u,v]}|v-u|^\be\vee N^{-\be}\leq 2^{-n\be}\vee N^{-\be}\|\del\Psi^X_N\|_{\be,N,[s,t]}
\]
if $[u,v]\in\cP_{s,t}^{(n)}$. Taking into account (\ref{2.10}) we obtain similarly to (\ref{3.2}) that
\begin{equation}\label{3.13}
|X_N(s,t)-\Psi_N^X(s,t)|\leq (1-2^{1-\be})^{-1}\|\del\Psi_N^X\|_{\be,N,[s,t]}(|t-s|^\be\vee N^{-\be})
\end{equation}
where, recall, $\be=3\al$.

Next, we set
\[
R^{\sig(X_N)}(s,u)=\sig(X_N)(s,t)-\nabla\sig(X_N(s))X_N(s,u)+\nabla\sig(X_N(s))R^{X_N}(s,u),
\]
\[
R^{X_N}(s,u)=X_N(s,u)-\sig(X_N(s))S_N(s,u)
\]
and define semi norms
\[
\|\cdot\|_{\al,N,h}=\sup_{I\subset[0,T],|I|\leq h}\|\cdot\|_{\al,N,I}\,\,\mbox{and}\,\,\|\cdot\|_{2\al,N,h}=\sup_{I\subset[0,T],|I|\leq h}\|\cdot\|_{2\al,N,I}.
\]
Proceeding in the same way as in Section \ref{subsec3.1} we obtain
\begin{eqnarray}\label{3.14}
&\|R^{X_N}\|_{2\al,N,h}\leq(1-2^{1-3\al})^{-1}(\frac 12\|\nabla^2\sig\|_\infty\| X_N\|^2_{\al,N,h}\| S_N\|_{\al,N,h}\\
&+(\|\nabla\sig\|_\infty^2+\|\nabla^2\sig\|_\infty\|\sig\|_\infty)\| X_N\|_{\al,N,h}\|\bbS_N\|_{2\al,N,h}+\| b\|_\infty (h\vee N^{-1})^{1-2\al}\nonumber\\
&+\|\nabla b\|_\infty\| X_N\|_{\al,N,h}(h\vee N^{-1})^{1-2\al})h^\al\vee N^{-\al}+\|\nabla\sig\|_\infty\|\sig\|_\infty\|\bbS_N\|_{2\al,N,h}.\nonumber
\end{eqnarray}

Set
\[
\hat C_\sig=4(\|\sig\|_\infty+\|\nabla^2\sig\|_\infty+\|\nabla\sig\|_\infty^2+\|\nabla^2\sig\|_\infty\|\sig\|_\infty+2\|\nabla\sig\|_\infty\|\sig\|_\infty+1).
\]
Taking into account the assumption (\ref{2.12}) we see that there exists $N_0\geq 1$ large enough such that for any $N\geq N_0$,
\[
N^{-\al}\hat C_\sig(\| S_N\|_{\al,N}+\sqrt{\|\bbS_N\|_{2\al,N}})+N^{-(1-2\al)}(\| b\|_\infty+\|\nabla b\|_\infty)\leq\frac 5{72}(1-2^{1-3\al}).
\]
Hence, for $N\geq N_0$,
\begin{eqnarray*}
&h_{0,N}=\sup\{ h\leq 1:\,\hat C_\sig(h\vee N^{-1})^{\al} (\| S_N\|_{\al,N}+\sqrt{\|\bbS_N\|_{2\al,N}})\\
&+(h\vee N^{-1})^{1-2\al}(\| b\|_\infty+\|\nabla b\|_\infty)\leq\frac 5{72}(1-2^{1-3\al})\}
\end{eqnarray*}
is well defined and $h_{0,N}\geq N^{-1}$.

As in Section \ref{subsec3.1}, for any $h\leq h_{0,N}$ we derive similarly to (\ref{3.9}) from (\ref{3.14}) and the equality
\[
X_N(s,t)=\sig(X_N(s))S_N(s,t)+R^{X_N}(s,t)
\]
that
\begin{eqnarray}\label{3.15}
&\,\,\quad\| X_N\|_{\al,N,h}\leq\|\sig\|_\infty\| S_N\|_{\al,N,h}+2(h^\al\vee N^{-\al})\| X_N\|^2_{\al,N,h}+2(h^\al\vee N^{-\al})\\
&+h^\al\vee N^{-\al}(2\|\nabla\sig\|_\infty\|\sig\|_\infty+1)\|\bbS_N\|_{2\al,N,h}+2(h^\al\vee N^{-\al})\| X_N\|_{\al,N,h}.
\nonumber\end{eqnarray}
Multiplying both parts of (\ref{3.15}) by $4(h^\al\vee N^{-\al})$ and setting
\begin{eqnarray*}
&\la_{h,N}=4(h^\al\vee N^{-\al})\|\sig\|_\infty\| S_N\|_{\al,N}\\
&+4(h^{2\al}\vee N^{-2\al})(2\|\nabla\sig\|_\infty\|\sig\|_\infty+1)\|\bbS_N\|_{2\al,N}+8(h^{2\al}\vee N^{-2\al})
\end{eqnarray*}
and $\vf_{h,N}=4(h^\al\vee N^{-\al})\|X_N\|_{\al,N,h}$ we obtain
\begin{equation}\label{3.16}
\vf_{h,N}\leq\la_{h,N}+\frac 12\vf^2_{h,N}+2(h^\al\vee N^{-\al})\vf_{h,N}.
\end{equation}

By the choice of $h_{0,N}$ above we see that $\la_{h,N}\leq \frac 5{72}$ for any $h\leq h_{0,N}$ and $N\geq N_0$. Since we also have that
$8(h^{2\al}\vee N^{-2\al})\leq\frac 5{72}$, it follows that $2(h^\al\vee N^{-\al})<1/4$ and (\ref{3.16}) implies that
\[
\vf_{h,N}\leq 2\la_{h,N}+\vf^2_{h,N}.
\]
In the same way as in Section \ref{subsec3.1} we conclude that for $h\leq h_{0,N}$ either $\vf_{h,N}\geq 5/6$ or $\vf_{h,N}\leq 1/6$.
But now we cannot argue that $\vf_{h,N}$ tends to zero as $h\to 0$ and by this reason we will estimate $\vf_{h,N}$ directly from the
definition.

Observe that when $0<t-s<N^{-1}$ then (\ref{2.11}) becomes $X_N(s,t)=\Psi^X_N(s,t)$, and so
\begin{equation}\label{3.17}
\vf_{h,N}\leq 4N^{-\al}\|\sig\|_\infty(\| S_N\|_{\al,N}+\|\nabla\sig\|_{\infty}N^{-\al}\|\bbS_N\|_{2\al,N}).
\end{equation}
This together with (\ref{2.12}) yields that there exists $N_1\geq N_0$ such that for all $N\geq N_1$ and $h<N^{-1}$ we have $\vf_{h,N}\leq 1/6$.
Arguing similarly to Section \ref{subsec3.1} we see that $\vf_h\leq 3\vf_{h/2}$ for all $h$ and by the monotonicity of $\vf_h$ in $h$
we conclude that it cannot jump from below $1/6$ to above $5/6$. Hence, $\vf_{h,N}\leq 1/6$ for all $h\leq h_{0,N}$  and $N\geq N_1$. It follows
similarly to Section \ref{subsec3.1} that
\begin{equation}\label{3.18}
\| X_N\|_{\al,N,h}\leq\tilde C_\sig(\| S_N\|_{\al,N}+\sqrt{\|\bbS_N\|_{2\al,N}}+1)
\end{equation}
provided $N\geq N_1$ and $h\leq h_{0,N}$, where $\tilde C_\sig$ is the same as in Section \ref{subsec3.1}.
Finally, we conclude similarly to Section \ref{subsec3.1} that (\ref{2.13}) holds true for all $N\geq N_1$ while for $N<N_1$ we obtain
(\ref{2.13}) by choosing a sufficiently large constant $C_T$ there.

\section{Completing proofs in deterministic setup}\label{sec4}\setcounter{equation}{0}
\subsection{Continuous time case}\label{subsec4.1}
Set $\Del(s,t)=\Psi^X(s,t)-\Psi^{\Xi}(s,t)$. Then by Assumption \ref{ass2.1},
\begin{equation}\label{4.1}
|X(s,t)-\Xi(s,t)-\sum_{i=1}^{m_n}\Del(t_{i-1},t_i)|\to 0\,\,\mbox{as}\,\, n\to\infty
\end{equation}
when $\cP_{s,t}^{(n)}=\{[t_{i-1}^{(n)},t_i^{(n)}],\, i=1,...,m_n\}$ is a sequence of partitions appearing
in Assumption \ref{ass2.1}. Let $\del\Del(s,u,t)=\Del(s,t)-\Del(s,u)-\Del(u,t)$. Then, in the same way as in
(\ref{3.2}) we obtain that
\begin{equation}\label{4.2}
|X(s,t)-\Xi(s,t)-\Del(s,t)|\leq (1-2^{1-\be})^{-1}\|\del\Del\|_{\be,[s,t]}|t-s|^\beta
\end{equation}
where $\be=3\al$ and, again,
\[
\|\del\Del\|_{\be,[s,t]}=\sup_{s\leq u<w<v\leq t}\frac {|\del\Del(u,w,v)|}{|v-u|^\be},\,\,\,[s,t]\subset[0,T].
\]
In the same way as in (\ref{3.3}) we have
\begin{eqnarray}\label{4.3}
&\del\Del(s,u,t)=-R^{\sig(X)}(s,u)S(u,t)+R^{\sig(\Xi)}(s,u)W(u,t)\\
&-\nabla\sig(X)\sig(X)(s,u)\bbS(u,t)+\nabla\sig(\Xi)\sig(\Xi)(s,u)\bbW(u,t)\nonumber\\
&-b(X)(s,u)(t-u)+b(\Xi)(s,u)(t-u)\nonumber
\end{eqnarray}
and similarly to (\ref{3.6+}) we conclude that $\|\del\Del\|_\be<\infty$.

Next, we have
\begin{eqnarray}\label{4.4}
&|\Del R(s,t)|=|R^X(s,t)-R^\Xi(s,t)|\leq |X(s,t)-\Xi(s,t)-\Del(s,t)|\\
&+|\nabla\sig(X(s))\sig(X(s))\bbS(s,t)-\nabla\sig(\Xi(s))\sig(\Xi(s))\bbW(s,t)|\nonumber\\
&+|b(X(s))-b(\Xi(s))|(t-s)\leq(1-2^{1-3\al})^{-1}\|\del\Del\|_{3\al,[s,t]}|t-s|^{3\al}\nonumber\\
&+(\|\nabla^2\sig\|_\infty\|\sig\|_\infty+\|\nabla\sig\|^2)|X(s)-\Xi(s)|\|\bbS\|_{2\al,[s,t]}|t-s|^{2\al}\nonumber\\
&+\|\nabla\sig\|_\infty\|\sig\|_\infty\|\bbS-\bbW\|_{2\al,[s,t]}|t-s|^{2\al}+\|\nabla b\|_\infty|X(s)-\Xi(s)|(t-s).\nonumber
\end{eqnarray}
By (\ref{4.3}),
\begin{eqnarray}\label{4.5}
&\|\del\Del\|_{2\al,[s,t]}\leq\|R^{\sig(X)}-R^{\sig(\Xi)}\|_{2\al,[s,t]}\| W\|_{\al,[s,t]}\\
&+\| R^{\sig(X)}\|_{2\al,[s,t]}\| S-W\|_{\al,[s,t]}+\|\nabla\sig(X)\sig(X)-\nabla\sig(\Xi)\sig(\Xi)\|_{\al.[s,t]}\|\bbS\|_{2\al,[s,t]}\nonumber\\
&+\|\nabla\sig(\Xi)\sig(\Xi)\|_{\al,[s,t]}\|\bbS-\bbW\|_{2\al,[s,t]}+\| b(X_N)-b(\Xi)\|_{\al,[s,t]}|t-s|^{1-2\al}.
\nonumber\end{eqnarray}
By the definition of $R^{\sig(X)}$ and $R^{\sig(\Xi)}$,
\begin{eqnarray}\label{4.6}
&\| R^{\sig(X)}-R^{\sig(\Xi)}\|_{2\al,[s,t]}\leq\|\sig(X)-\nabla\sig(X)X-\sig(\Xi)+\nabla\sig(\Xi)\Xi\|_{2\al,[s,t]}\\
&+\|\nabla^2\sig\|_\infty|X(s)-\Xi(s)|\| R^X\|_{2\al,[s,t]}+\|\nabla\sig\|_\infty\|\Del R\|_{2\al,[s,t]}.\nonumber
\end{eqnarray}

Employing the Taylor formula with two terms and the integral remainder we obtain that for $s\leq u<v\leq t$,
\begin{eqnarray}\label{4.7}
&|b(X_N)(s,u)-b(\Xi)(s,u)|=|b(X(u))-b(X(s))-(b(\Xi(u))-b(\Xi(s)))|\\
&\leq|\int_0^1(\nabla b(X(s)+\te X(s,u))X(s,u)-\nabla b(\Xi(s)+\te\Xi(s,u))\Xi(s,u))d\te|\nonumber\\
&\leq\|\nabla^2b\|_\infty |X(s,u)||X(s)-\Xi(s)|+(\|\nabla^2b\|_\infty|X(s,u)|\nonumber\\
&+\|\nabla b\|_\infty)|X(s,u)-\Xi(s,u)|\nonumber
\end{eqnarray}
and
\begin{eqnarray}\label{4.8}
&|\nabla\sig(X)\sig(X)(u,v)-\nabla\sig(\Xi)\sig(\Xi)(u,v)|\\
&\leq|\int_0^1(\nabla(\nabla\sig(X)\sig(X))(X(u)+\te X(u,v))X(u,v)\nonumber\\
&-\nabla(\nabla\sig(\Xi)\sig(\Xi))(\Xi(u)+\te\Xi(u,v)\Xi(u,v))d\te|\leq(\|\nabla^3\sig\|_\infty\|\sig\|_\infty\nonumber\\
&+\|\nabla^2\sig\|_\infty\|\nabla\sig\|_\infty+2\|\nabla^2\sig\|_\infty)(|X(u)-\Xi(u)|+|X(u,v)-\Xi(u,v)|)|X(u,v)|\nonumber\\
&(\|\nabla^2\sig\|_\infty\|\sig\|_\infty+\|\nabla\sig\|^2_\infty)|X(u,v)-\Xi(u,v)|.\nonumber
\end{eqnarray}
Similarly, by the Taylor formula with three terms and the integral remainder we obtain
\begin{eqnarray}\label{4.9}
&|\sig(X)(u,v)-\nabla\sig(X(u))X(u,v)-\sig(\Xi)(u,v)+\nabla\sig(\Xi(u))\Xi(u,v)|\\
&=\frac 12|\int_0^1(\nabla^2\sig(X(u)+\te X(u,v))X(u,v)\otimes X(u,v)\nonumber\\
&-\nabla^2\sig(\Xi(u)+\te\Xi(u,v))\Xi(u,v)\otimes\Xi(u,v))d\te|\nonumber\\
&\leq\frac 12\|\nabla^3\sig\|_\infty(|X(u)-\Xi(u)|+|X(u,v)-\Xi(u,v)|)|X(u,v)|^2\nonumber\\
&+\|\nabla^2\sig\|_\infty(|X(u,v)|+|\Xi(u,v)|)|X(u,v)-\Xi(u,v)|).
\nonumber\end{eqnarray}

Combining (\ref{4.4})--(\ref{4.9}) we obtain
\begin{eqnarray}\label{4.10}
&\|\Del R\|_{2\al,[s,t]}\leq L_1(s,t)\|\Del R\|_{2\al,[s,t]}+L_2(s,t)\| X-\Xi\|_{\al,[s,t]}\\
&+(L_3(s,t)+L_4\|\bbS\||_{2\al})|X(s)-\Xi(s)|+L_5(s,t)+\|\nabla\sig\|_\infty\|\sig\|_\infty\|\bbS-\bbW\|_{2\al}\nonumber
\end{eqnarray}
where
\[
L_1(s,t)=(1-2^{1-3\al})^{-1}|t-s|^\al\|\nabla\sig\|_\infty\| W\|_\al,
\]
\begin{eqnarray*}
&L_2(s,t)=(1-2^{1-3\al})^{-1}|t-s|^\al\big(\frac 12|t-s|^\al\|\nabla^3\sig\|_\infty\| X\|^2_{\al,[s,t]}\| W\|_\al\\
&+\|\nabla^2\sig\|_\infty(\| X\|_{\al,[s,t]}+\|\Xi\|_{\al,[s,t]})\| W\|_\al+|t-s|^\al(\|\nabla^3\sig\|_\infty\|\sig\|_\infty\\
&+\|\nabla^2\sig\|_\infty\|\nabla\sig\|_\infty+2\|\nabla^2\sig\|_\infty)\| X\|_{\al,[s,t]}\|\bbS\|_{2\al}+(\|\nabla^2\sig\|_\infty\|\sig\|_\infty\\
&+\|\nabla\sig\|^2_\infty)\|\bbS\|_{2\al}+(\|\nabla^2b\|\| X\|_{\al,[s,t]}|t-s|^{1-\al}+
\|\nabla b\|_\infty|t-s|^{1-2\al}\big),
\end{eqnarray*}
\begin{eqnarray*}
&L_3(s,t)=(1-2^{1-3\al})^{-1}|t-s|^\al\big(\frac 12\|\nabla^3\sig\|_\infty\| X\|^2_{\al,[s,t]}\| W\|_\al\\
&+\|\nabla^2\sig\|\| R^X\|_{2\al,[s,t]}\| W\|_\al+(\|\nabla^3\sig\|_\infty\|\sig\|_\infty+\|\nabla^2\sig\|_\infty\|\nabla\sig\|_\infty\\
&+2\|\nabla^2\sig\|_\infty)\| X\|_{\al,[s,t]}\|\bbS\|_{2\al}+\|\nabla^2b\|_\infty\| X\|_{\al,[s,t]}|t-s|^{1-2\al}\big)+\|\nabla b\|_\infty|t-s|^{1-2\al},\\ &L_4=\|\nabla^2\sig\|_\infty\|\sig\|_\infty+\|\nabla\sig\|^2_\infty\quad\mbox{and}\\
&L_5(s,t)=(1-2^{1-3\al})^{-1}|t-s|^\al\big(\| R^{\sig(X)}\|_{2\al,[s,t]}\| S-W\|_\al\\
&+(\|\nabla^2\sig\|_\infty\|\sig\|_\infty+\|\nabla\sig\|^2_\infty)\|\Xi\|_\al\|\bbS-\bbW\|_{2\al,[s,t]}\big).
\end{eqnarray*}

On the other hand,
\[
X(s,t)-\Xi(s,t)=(\sig(X(s))-\sig(\Xi(s)))S(s,t)+\sig(\Xi(s))(S(s,t)-W(s,t))+\Del R(s,t),
\]
and so
\begin{equation}\label{4.11}
\| X-\Xi\|_{\al,[s,t]}\leq\|\nabla\sig\|_\infty|X(s)-\Xi(s)|\| S\|_\al+\|\sig\|_\infty\| S-W\|_\al+\|\Del R\|_{2\al,[s,t]}|t-s|^\al.
\end{equation}

Taking into account (\ref{3.6}), (\ref{3.8}), (\ref{3.12}), (\ref{4.10}) and (\ref{4.11}) we conclude that there exists a small constant $\rho_{\al,\sig}>0$
depending only on $\al$ and $\sig$ but not on $S,\bbS,X,\Xi, W,\bbW$ (and which can be estimated from (\ref{3.6}), (\ref{3.8}), (\ref{3.12}) and
(\ref{4.11})) such that if $|t-s|$ is small enough and, in particular,
\begin{equation*}
|t-s|^\al(\| S\|_\al+\sqrt{\|\bbS\|_{2\al}}+\| W\|_\al+\sqrt{\|\bbW\|_{2\al}}+1)\leq\rho_{\al,\sig}
\end{equation*}
then
\begin{eqnarray*}
&L_1(s,t)\leq 1/2,\, L_2(s,t)(t-s)^\al\leq 1/4,\, (L_3(s,t)+L_4\|\bbS\|_{2\al})(t-s)^\al\\
&\leq\| S\|_\al+\sqrt{\|\bbS\|_{2\al}}+\| W\|_\al+\sqrt{\|\bbW\|_{2\al}},\quad\mbox{and}\\
&(L_5(s,t)+\|\nabla\sig\|_\infty\|\sig\|_\infty\|\bbS-\bbW\|_{2\al})(t-s)^\al\leq\| S-W\|_\al+\sqrt{\|\bbS-\bbW\|_{2\al}}.
\end{eqnarray*}
Combining this with (\ref{4.9}) we obtain that
\begin{eqnarray*}
&\|\Del R\|_{2\al,[s,t]}|t-s|^\al\leq\frac 12\| X-\Xi\|_{\al,[s,t]}+(\| S\|_\al+\sqrt{\|\bbS\|_{2\al}}\\
&+\| W\|_\al+\sqrt{\|\bbW\|_{2\al}})|X(s)-\Xi(s)|+\| S-W\|_\al+\sqrt{\|\bbS-\bbW\|_{2\al}}.
\end{eqnarray*}
This together with (\ref{4.10}) yields
\begin{eqnarray}\label{4.12}
&\| X-\Xi\|_{\al,[s,t]}\leq(1+\|\nabla\sig\|_\infty)(\| S\|_\al+\sqrt{\|\bbS\|_{2\al}}\\
&+\| W\|_\al+\sqrt{\|\bbW\|_{2\al}})|X(s)-\Xi(s)|+\| S-W\|_\al+\sqrt{\|\bbS-\bbW\|_{2\al}}.\nonumber
\end{eqnarray}

Set
\[
h_2=\max\{ h\leq h_1:\, h^\al(\| S\|_\al+\sqrt{\|\bbS\|_{2\al}}+\| W\|_\al+\sqrt{\bbW\|_{2\al}}+1)\leq\rho_{\al,\sig}\},
\]
where $h_1$ was defined in Section \ref{subsec3.1},
\begin{eqnarray*}
&A=(1+\|\nabla\sig\|_\infty)(\| S\|_\al+\sqrt{\|\bbS\|_{2\al}}+\| W\|_\al+\sqrt{\|\bbW\|_{2\al}})\\
&\mbox{and}\quad B=\| S-W\|_\al+\sqrt{\|\bbS-\bbW\|_{2\al}}.
\end{eqnarray*}
By our definition of the $\al$-H\" older (semi) norms
\begin{eqnarray*}
&\| X-\Xi\|_{\al,[0,s]}=\sup_{0\leq u<v\leq s}\frac {(X(u,v)-\Xi(u,v)|}{|v-u|^\al}\geq\frac {|X(0,s)-\Xi(0,s)|}{s^\al}\\
&\geq T^{-\al}(|X(s)-\Xi(s)|-|X(0)-\Xi(0)|),
\end{eqnarray*}
and so
\[
|X(s)-\Xi(s)|\leq T^\al\| X-\Xi\|_{\al,[0,s]}+|X(0)-\Xi(0)|.
\]
Also, it follows easily from the definition of our (semi) norms that for $0\leq s<t\leq T$,
\[
\| X-\Xi\|_{\al,[0,t]}\leq\| X-\Xi\|_{\al,[0,s]}+\| X-\Xi\|_{\al,[s,t]}.
\]
Then we have
\[
\| X-\Xi\|_{\al,[0,(n+1)h_2]}\leq (AT^\al+1)\| X-\Xi\|_{\al,[0,nh_2]}+A|X(0)-\Xi(0)|+B
\]
and (\ref{2.7}) follows by induction taking into account that the number of induction steps $[T/h_2]+1$ has the order of
\[
(\| S\|_\al+\sqrt{\|\bbS\|_{2\al}}+\| W\|_\al+\sqrt{\|\bbW\|_{2\al}}+1)^{1/\al}.
\]

\subsection{Discrete time case}\label{subsec4.2}
In fact, we have here a "mixed" time case as $\Psi^X_N$ comes from the discrete time construction while $\Psi^\Xi$ comes
from the continuous time one. So, set $\Del_N(s,t)=\Psi^X_N(s,t)-\Psi^\Xi(s,t)$. Then by the assumptions (\ref{2.5}) and (\ref{2.10}),
\begin{equation}\label{4.13}
|X_N(s,t)-\Xi(s,t)-\sum_{i=1}^{m_n}\Del_N(t_{i-1},t_i)|\to 0\,\,\mbox{as}\,\, n\to\infty
\end{equation}
when $\cP_{s,t}^{(n)}=\{[t^{(n)}_{i-1},t_i^{(n)}],\, i=1,...,m_n\}$ is a sequence of partitions constructed before Assumption \ref{ass2.1}.
Let $\del\Del_N(s,u,t)=\Del_N(s,t)-\Del_N(s,u)-\Del_N(u,t)$. Then, similarly to (\ref{3.13}) we obtain
\begin{equation}\label{4.14}
|X_N(s,t)-\Xi(s,t)-\Del_N(s,t)|\leq (1-2^{1-\be})^{-1}\|\del\Del_N\|_{\be,N,[s,t]}(|t-s|^\be\vee N^{-\be})
\end{equation}
where, as before, $\be=3\al$.

In the same way as in (\ref{3.3}) we have
\begin{eqnarray}\label{4.15}
&\del\Del_N(s,u,t)=-R^{\sig(X_N)}(s,u)S_N(u,t)+R^{\sig(\Xi)}(s,u)W(u,t)\\
&-\nabla\sig(X_N)\sig(X_N)(s,u)\bbS_N(u,t)+\nabla\sig(\Xi)\sig(\Xi)(s,u)\bbW(u,t)\nonumber\\
&-b(X_N)(s,u)N^{-1}([Nt]-[Nu])+b(\Xi)(s,u)(t-u).\nonumber
\end{eqnarray}
Proceeding similarly to (\ref{4.4}) but taking into account that (\ref{4.14}) and (\ref{4.15}) have slightly different
from (\ref{4.2}) and (\ref{4.3}) form, we obtain that
\begin{eqnarray}\label{4.16}
&|\Del_NR(s,t)|=|R^{X_N}(s,t)-R^{\Xi}(s,t)|\\
&\leq(1-2^{1-3\al})^{-1}\|\del\Del_N\|_{3\al,N,[s,t]}(|t-s|^{3\al}\vee N^{-3\al})\nonumber\\
&+\|\nabla^2\sig\sig+(\nabla\sig)^2\|_\infty|X_N(s)-\Xi(s)|\|\bbS_N\||_{2\al,N,[s,t]}(|t-s|^{2\al}\vee N^{-2\al})\nonumber\\
&+\|\nabla\sig\sig\|_\infty\|\bbS_N-\bbW\|_{2\al,N,[s,t]}(|t-s|^{2\al}\vee N^{-2\al})\nonumber\\
&+\|\nabla b\|_\infty|X_N(s)-\Xi(s)||t-s|+2\| b\|_\infty N^{-1}\nonumber
\end{eqnarray}
where we use that
\[
|N^{-1}([Nt]-[Ns])-(t-s)|\leq 2N^{-1}.
\]

Next, we continue similarly to (\ref{4.5})--(\ref{4.8}) arriving at the analogy of the inequality (\ref{4.9}) having now the form
\begin{eqnarray}\label{4.17}
&\|\Del_NR\|_{2\al,N,[s,t]}\leq L_{1,N}(s,t)+L_{2,N}(s,t)\| X_N-\Xi\|_{\al,N,[s,t]}\\
&+(L_{3,N}(s,t)+L_{4,N}\|\bbS_N\|_{2\al,N})|X_N(s)-\Xi(s)|+L_{5,N}(s,t)\nonumber\\
&+\|\sig\|_\infty\|\sig\|_\infty\|\bbS_N-\bbW\|_{2\al,N}+2\| b\|_\infty N^{-(1-2\al)}
\nonumber\end{eqnarray}
where $L_{i,N}$ has the same form as $L_i,\, i=1,...,5$ in (\ref{4.9}) with $X,\, S$ and $\bbS$ replaced by $X_N,\, S_N$ and $\bbS_N$,
respectively. The estimate (\ref{4.11}) becomes here
\begin{eqnarray}\label{4.18}
&\| X_N-\Xi\|_{\al,N,[s,t]}\leq\|\nabla\sig\|_\infty|X_N(s)-\Xi(s)|\| S\|_{\al,N}\\
&+\|\sig\|_\infty\| S_N-W\|_{\al,N}+\|\Del R\|_{2\al,N,[s,t]}(|t-s|^\al\vee N^{-1}).\nonumber
\end{eqnarray}

Taking into account (\ref{3.14}), (\ref{3.18}), (\ref{4.17}) and (\ref{4.18}) we see as in Section \ref{subsec4.1}
 that there exists a small constant $\rho_{\al,\sig}>0$, depending only on $\al$ and $\sig$, such that if $|t-s|\vee N^{-1}$
 is small enough and, in particular,
\begin{equation}\label{4.19}
(|t-s|^\al\vee N^{-\al})(\| S_N\|_{\al,N}+\sqrt{\|\bbS_N\|_{2\al,N}}+\| W\|_\al+\sqrt{\|\bbW\|_{2\al}}+1)\leq\rho_{\al,\sig}
\end{equation}
then
\begin{eqnarray*}
&L_{1,N}(s,t)\leq 1/2,\, L_{2,N}(s,t)(|t-s|^\al\vee N^{-\al})\leq 1/4,\\
& (L_{3,N}(s,t)+L_{4,N}\|\bbS_N\|_{2\al,N})(|t-s|^\al\vee N^{-\al})\leq\| S_N\|_{\al,N}+\sqrt{\|\bbS_N\|_{2\al,N}}\\
&+\| W\|_\al+\sqrt{\|\bbW\|_{2\al}},\quad\mbox{and}\quad(L_{5,N}(s,t)+\|\nabla\sig\|_\infty\|\sig\|_\infty\|\bbS_N-\bbW\|_{2\al,N})\\
&\times(|t-s|^\al\vee N^{-\al})\leq\| S_N-W\|_{\al,N}+\sqrt{\|\bbS_N-\bbW\|_{2\al,N}}.
\end{eqnarray*}
Now observe that unlike in Section \ref{subsec4.1} we cannot have (\ref{4.19}) just by choosing $|t-s|$ small enough since we have, in fact,
$|t-s|\vee N^{-1}$ instead, but in view of the assumption (\ref{2.12}) we can find $N_2\geq N_1$ large enough so that (\ref{4.19}) holds true
for all $N\geq N_2$ and $|t-s|\geq N^{-1}$ small enough.

It is clear that it suffices to prove (\ref{2.14}) only for all $N\geq N_2$ since for $N<N_2$ we can obtain (\ref{2.14}) just by adjusting
a constant $C_T$ there. Thus we proceed assuming that (\ref{4.19}) holds true and we arrive at the estimate similar to (\ref{4.12}),
\begin{eqnarray}\label{4.20}
&\| X_N-\Xi\|_{\al,N,[s,t]}\leq(1+\|\nabla\sig\|_\infty)(\| S_N\|_{\al,N}+\sqrt{\|\bbS_N\|_{2\al,N}}\\
&+\| W\|_\al+\sqrt{\|\bbW\|_{2\al}})|X_N(s)-\Xi(s)|+\| S_N-W\|_{\al,N}+\sqrt{\|\bbS_N-\bbW\|_{2\al,N}}.\nonumber
\end{eqnarray}
We conclude the proof of (\ref{2.14}) by induction in the same way as in Section \ref{subsec4.1}.

Next, we will establish Corollary \ref{cor2.4}. For any partition $\cP=\{[t_{i-1},t_i],\, i=1,...,m,\, 0=t_0<t_1<...<t_m=T\}$ of the
interval $[0,T]$ we can write
\[
\sum_{i=1}^m|X_N(t_{i-1},t_i)-\Xi(t_{i-1},t_i)|^p\leq J_1+2^{p-1}(J_2+J_3)
\]
where taking into account that $\al p>1$,
\begin{eqnarray*}
& J_1=\sum_{i:|t_i-t_{i-1}|\geq 1/N}|X_N(t_{i-1},t_i)-\Xi(t_{i-1},t_i)|^p\\
&\leq\| X_N-\Xi\|^p_{\al,N}\sum_{i=1}^m|t_i-t_{i-1}|^{\al p}\leq T^{\al p}\| X_N-\Xi\|^p_{\al,N},\\
& J_2=\sum_{i:0<|t_i-t_{i-1}|< 1/N}|X_N(t_{i-1},t_i)|^p\\
&\leq\sum_{k=1}^{[TN]}|X_N(\frac {k-1}N,\frac kN)|^p\leq TN^{-(\al p-1)}\|X_N\|_{\al,N}\quad\mbox{and}\\
&J_3=\sum_{i:0<|t_i-t_{i-1}|< 1/N}|\Xi(t_{i-1},t_i)|^p\\
&\leq\|\Xi\|_\al\sum_{i:0<|t_i-t_{i-1}|< 1/N}|t_i-t_{i-1}|^{\al p}\leq TN^{-(\al p-1)}\|\Xi\|_\al.
\end{eqnarray*}
These together with
\[
(J_1+2^{p-1}(J_2+J_3))^p\leq J_1^{1/p}+2(J_2^{1/p}+J_3^{1/p})
\]
yields (\ref{2.15}).

\section{Probabilistic setup  }\label{sec5}\setcounter{equation}{0}
\subsection{Discrete time case}\label{subsec5.1}
First, we will check the conditions of Theorem \ref{thm2.3} in the probabilistic setup when $X_N(t)$ is given by the sum (\ref{2.16}) built with a stationary sequence $\xi$ while $S_N$ and $\bbS_N$ are the sums and iterated sums made of 
this sequence. Because of the definition of $\|\cdot\|_{\be,N}$ the condition (\ref{2.8}) is clearly satistied.
It is immediate to check that sums and the corresponding iterated sums satisfy the Chen relation (\ref{2.9}). Of course, 
the condition (\ref{2.10}) is clear
 and, in fact, we have here (\ref{2.11}). Indeed, if $|u-v|<1/N$ then $S_N(u,v)=0$ if $\frac kN\leq u\leq v<\frac {k+1}N$ for some integer $k$ and
 $S_N(u,v)=\xi(k)$ if $\frac kN\leq u<\frac {k+1}N\leq v$ while $\bbS_N(u,v)=0$ in both cases. Hence, if $|v-u|<\frac 1N$ then $\Psi^X_N(u,v)=0$
 in the first case and $\Psi_N^X(u,v)=\sig(X_N(u))\xi(k)$ in the second case which gives (\ref{2.11}). This implies 
 also that (\ref{2.8}) is satisfied for $R_N(s,t)=X_N(s,t)-\sig(X_N(s))S_N(s,t)$.

Next, the diffusion $\Xi=\Xi_N$ solving the stochastic differential equation (\ref{2.18}) can be taken to satisfy (\ref{2.5})  almost surely
since stochastic integrals can be obtained as limits of integral sums (for instance, in the sense of rough integration, see Chapter 4 in \cite{FH}).
The a.s. $\al$-H\" older continuity of $W_N$ and $2\al$-H\" older continuity of $\bbW_N$ for $\al\in(\frac 13,\frac 12)$, required in (\ref{2.2}),
follows from the Kolmogorov--Chentsov type theorem (see Theorem 3.1 in \cite{FH} and its generalization Theorem 3.3 in \cite{Ki23}). We can use the same theorem for $\Xi=\Xi_N$ estimating its moments (see \cite{Mao}),
\begin{equation*}
E|\Xi_N(s,t)|^{2M}\leq 2^{2M}\big(C(M(2M-1))^M\|\sig\|_\infty^{2M}(t-s)^M+(\| b\|_\infty+\| c\|_\infty)^{2M}(t-s)^{2M}\big)
\end{equation*}
where a constant $C>0$ comes here since we take the stochastic integral with respect to a Brownian motion with some
covariance matrix. By the Kolmogorov--Chentsov type theorem mentioned above this estimate
 implies existence of an $\al$-H\" older continuous modification of $\Xi$ for each $\al<1/2$ taking $M>(\frac 12-\al)^{-1}$.

We have to check also that
\begin{equation}\label{5.1}
\| R^{\Xi_N}\|_{2\al}<\infty\quad\mbox{a.s.}
\end{equation}
where
\[
R^{\Xi_N}(s,t)=\Xi_N(s,t)-\sig(\Xi_N(s))W_N(s,t)=R_1^{\Xi_N}(s,t)+R_2^{\Xi_N}(s,t)
\]
with
\[
R_1^{\Xi_N}(s,t)=\int_s^t\big(\sig(\Xi_N(u))-\sig(\Xi_N(s))\big)dW_N
\]
and
\[ 
 R_2^{\Xi_N}(s,t)=\int_s^t\big(b(\Xi_N(u))+c(\Xi_N(u))\big)du.
\]
Since $\al\in(\frac 13,\,\frac 12)$ it is clear that $\| R_2^{\Xi_N}\|_{2\al}<\infty$.

Concerning $R_1^{\Xi_N}$ we have by the standard inequalities for stochastic integrals ( see, for instance,
Theorem 7.1 in \cite{Mao}) that for all $s,t\in[0,T],\, s<t$ and $M\geq 1$,
\begin{eqnarray*}
&E|R_1^{\Xi_N}(s,t)|^{2M}\leq C_{N,T,M}(t-s)^{M-1}\int_s^tE|\Xi_N(u)-\Xi_N(s)|^{2M}du\\
&\leq\tilde C_{N,T,M}(t-s)^{M(2\al+1)}E\|\Xi_N\|^{2M}_\al
\end{eqnarray*}
where $C_{N,T,M},\tilde C_{N,T,M}>0$ are some constants. Now, by the Kolmogorov--Chentsov type theorem 
(see Theorem 3.1 in \cite{FH} and Theorem 3.3 in \cite{Ki23}) we obtain that there exists an $\al$--H\" older 
continuous modification of $\Xi_N$ preserving its distributions and such that 
$R_1^{\Xi_N}$ is $\al+\frac 12-\frac 1{2M}$--H\" older continuous. Choosing $M>(1-2\al)^{-1}$ we see that (\ref{5.1}) holds true.

Next, we could obtain the condition (\ref{2.12}) from the law of iterated logarithm in the H\" older norm which will be discussed later on but for
now we will obtain (\ref{2.12}) by a more direct argument. By Lemma 3.2 from \cite{Ki23} it follows that for any $0\leq s\leq t\leq T$ with $t-s>1/N$,
\begin{equation}\label{5.2}
E|S_N(s,t)|^{2M}\leq C(M)(t-s)^M\,\,\,\mbox{and}\,\,\, E|\bbS_N(s,t)|^{2M}\leq C(M)(t-s)^{2M}
\end{equation}
where $C(M)>0$ does not depend on $N$. By the discrete time version of the Kolmogorov--Chentsov type theorem appearing as Proposition 3.4 in \cite{Ki23}
it follows from (\ref{5.2}) that
\begin{equation}\label{5.3}
E\| S_N\|^{2M}_{\al,N}\leq\tilde C_\al(M)<\infty\,\,\,\mbox{and}\,\,\, E\|\bbS_N\|^M_{2\al,N}\leq\tilde C_\al(M)<\infty
\end{equation}
for any $\al<1/2$, where $\tilde C_\al(M)>0$ does not depend on $N$. Hence, by the Chebyshev inequality
\[
P\{\| S_N\|_{\al,N}>N^{\al/2}\}\leq\tilde C_\al(M)N^{-\al M}
\]
and
\[
P\{\|\bbS_N\|_{2\al,N}>N^{\al}\}\leq\tilde C_\al(M)N^{-\al M}.
\]
Taking $\al>1/3$ and $M>3$ we see that these probabilities are estimated by converging in $N$ sequences, and so by the Borel--Cantelli lemma we
conclude that
\[
\| S_N\|_{\al,N}+\sqrt{\|\bbS_N\|_{2\al,N}}=O(N^{\al/2})\quad\mbox{a.s.}
\]
implying (\ref{2.12}) with probability one.

The conditions of Theorem \ref{thm2.3} are now satisfied and we obtain the estimates (\ref{2.13})--(\ref{2.15}) with probability one. In order to
derive Theorem \ref{thm2.5} we have to estimate modified H\" older (semi) norms appearing in (\ref{2.13})--(\ref{2.15}) taking into account that
in \cite{Ki23+} the corresponding estimates are done in the $p$-variation norm. It was shown in Theorem 2.1 from \cite{Ki23+} that the process $\xi$
can be redefined preserving its distribution on a sufficiently rich probability space which contains also a $d$-dimensional Brownian motion $\cW$
with the covariance matrix $\vs$ such that the rescaled process $W_N(t)=N^{-1/2}\cW(Nt)$ and the corresponding iterated integral process $\bbW_N$,
 defined in Section \ref{subsec2.3}, together with the redefined $S_N$ and $\bbS_N$ satisfy
 \begin{equation}\label{5.4}
 \| S_N-W_N\|_p=O(N^{-\ve})\,\,\,\mbox{and}\,\,\,\|\bbS_N-\bbW_N\|_{p/2}=O(N^{-\ve})\,\,\,\mbox{a.s.}
 \end{equation}
 for some $\ve>0$.

 In order to obtain similar estimates in the modified H\" older semi norm we will use the following observation. Clearly, for any map $V$ from
 $\{ 0\leq u\leq v\leq T\}$ to $\bbR^n$ or to $\bbR^n\otimes\bbR^n$,
\begin{eqnarray*}
&\| V\|_p^p\geq\sup_{\cP=\{[t_{i-1},t_i],\,|t_i-t_{i-1}|\geq N^{-\gam},\. i=1,...,m\}}\sum_{i=1}^m|V(t_{i-1},t_i)|^p\\
&\geq N^{-\gam\al p}\sup_{|t-s|\geq N^{-\gam}}\frac {|V(s,t)|^p}{|t-s|^{\al p}},\\
\end{eqnarray*}
where we take only one term in the sum but with the maximal ratio in the right hand side, and so
\begin{equation*}
\| V\|_p\geq N^{-\gam\al}\sup_{|t-s|\geq N^{-\gam}}\frac {|V(s,t)|}{|t-s|^{\al }}.
\end{equation*}

 Hence, the $\al$-H\" older norm of $V$ can be estimated as
 \[
 \| V\|_{\al,N}\leq\max(\| V\|_{\al,N,N^{-\gam}},\, N^{\gam\al}\| V\|_p)
 \]
 where, recall, $\|\cdot\|_{\al,N,h}$ is the semi norm introduced in Section \ref{subsec3.2} where we take the supremum over ratios
 $\frac {V(s,t)}{|t-s|^\al}$ where $N^{-1}\leq |t-s|\leq h$. Let $1/3<\al<\be<1/2$ and $0<h\leq 1$, then
 \[
 \| V\|_{\al,N,h}\leq \| V\|_{\be,N,h}h^{\be-\al}\leq \| V\|_{\be,N},
 \]
 and so
 \begin{equation}\label{5.5}
 \| V\|_{\al,N}\leq\max( N^{-(\be-\al)\gam}\| V\|_{\be,N},\,N^{\gam\al}\| V\|_p).
 \end{equation}

 Now, we apply (\ref{5.5}) to $V=S_N-W_N$  obtaining
 \begin{equation}\label{5.6}
 \| S_N-W_N\|_{\al,N}\leq\max(N^{-(\be-\al)\gam}(\| S_N\|_{\be,N}+\| W_N\|_{\be,N}),\, N^{\gam\al}\| S_N-W_N\|_p)
 \end{equation}
 and, similarly,
\begin{equation}\label{5.7}
 \|\bbS_N-\bbW_N\|_{2\al,N}\leq\max(N^{-4(\be-\al)\gam}(\| \bbS_N\|_{2\be,N}+\| \bbW_N\|_{2\be,N}),\, N^{2\gam\al}\| \bbS_N-\bbW_N\|_{p/2}).
 \end{equation}
 Choosing $\gam>0$ so small that $2\gam\al<\ve$ we will use (\ref{5.4}) to make the second terms of the maximums in (\ref{5.6}) and (\ref{5.7})  of order $O(N^{-\del})$ a.s. for some $\del>0$. The same estimate holds true for the 
 first terms there which can be obtained in the same way as the proof
 of (\ref{2.12}) above or by applying the law of iterated logarithm for $\| S_N\|_{\be,N}$, $\| W_N\|_{\be,N}$, $\|\bbS_N\|_{2\be,N}$ and $\|\bbW_N\|_{2\be,N}$ which will be discussed below.

 In order to complete the proof of Theorem \ref{thm2.5} it remains to show that the local Lipschitz coefficient in (\ref{2.14}) in front of
 $\| S_N-W\|_{\al,N}+\sqrt{\|\bbS_N-\bbW_N\|_{2\al,N}}$ can only grow in $N$ much slower than any power of $N$. This will follow from a version of the law of iterated logarithm in H\" older norms and since we already know that
 \[
 \| S_N-W_N\|_{\al,N}=O(N^{-\del})\quad\mbox{and}\quad\|\bbS_N-\bbW_N\|_{2\al,N}=O(N^{-\del})\,\,\mbox{a.s.},
 \]
 it suffices to establish this law only for $W_N$ and $\bbW_N$. In \cite{BBK} it was shown that all limit points 
 in H\" older norms as $N\to\infty$ of the vector functions
 \[
 \zeta_N(t)=\frac {W_N(t)}{\sqrt{\log\log N}}=\frac {\cW(Nt)}{\sqrt{N\log\log N}},\,\, t\in[0,T]
 \]
 belong to a compact set in $\|\cdot\|_\al$-norm topology for any $\al\in(1/3,1/2)$. In \cite{Bal} the corresponding
 law of iterated logarithm in the supremum norm was established for all iterated stochastic integrals of the form
 \[
 \bbW^{(n)}(t)=\int_{0\leq t_1\leq t_2\leq...\leq t_{n-1}\leq t}dW(t_1)\otimes dW(t_2)\otimes\cdots\otimes dW(t_{n})
 \]
 (see also \cite{Ki23+}).
 For our purposes it suffices to extend such law of iterated logarithm to H\" older norms just for $n=2$ which can
 be done similarly to \cite{BBK} relying on some large deviations arguments. It seems that it would be useful to obtain
 the law of iterated logarithm in H\" older norms for all $\bbW^{(n)},\, n\geq 2$ which we leave for another paper.
 Representing $\bbW^{(n)},\, n\geq 2$ as a diffusion process similarly to \cite{Bal} this should be possible to do extending
  the arguments in \cite{BBK} from Brownian motions to all diffusions (cf. \cite{Bal}) relying, in particular, on large
  deviations in H\" older norms for diffusions (see Theorem 19.9 in \cite{FV}).
 Observe that in the uniform norm the law of iterated logarithm for iterated sums and integrals was obtained in \cite{Ki23+}.
 
 In fact, for our purposes here we do not need the precise law of iterated logarithm and an estimate on $[0, T]$ of the
 form 
 \begin{equation}\label{5.8}
 \|\bbW_N\|_{2\al,N}\leq C(T)\log\log N,
 \end{equation}
 would suffice, where $C(T)>0$ is a random variable which depends on $T$ and $\al$ but not on $N$. This was proved
 in Proposition 6.3 from \cite{FK} for the $p$-variation norm. The above inequality in the H\" older norm can be proved similarly using the Strassen theorem in H\" older norms from \cite{BBK}.


 \subsection{Continuous time case}\label{subsec5.2}
 The condition (\ref{2.4}) of Assumption \ref{ass2.1} follows immediately from the properties of the integral, say, in the Stieltyes sense, as well as
 in the rough paths sense. The almost sure $\al$-H\" older continuity of $X_N$ follows directly from its definition as a solution of (\ref{2.21}). The almost sure finiteness of $\al$-H\" older norm of $S_N$ and $2\al$-H\" older norm of 
 $\bbS_N$ follows from the Kolmogorov--Chentsov theorem (see Theorem 3.1 in \cite{FH} and
 Theorem 3.3 in \cite{Ki23}) together with the corresponding moment bounds from Lemma 3.2 in \cite{Ki23}.

 By Theorems 2.2 and 2.3 from \cite{Ki23+} we have (\ref{5.4}) for both continuous time setups. In the same way as in Section \ref{subsec5.1} we extend these estimates to the $\al$-H\" older norms $\|\cdot\|_{\al}$ in the continuous time case. 
 As explained in Section \ref{subsec5.1} an extension of the law of iterated logarithm for $W_N$ in the $\al$-H\" older 
 norms from \cite{BBK} to $\bbW_N$ or an estimate of the form (\ref{5.8}) will
  show that the coefficient in front of $A_N=\| S_N-W_N\|_\al+\sqrt{\|\bbS_N-\bbW_N\|_{2\al}}$ (all items in (\ref{2.7}) should be supplied by the index $N$
   now except for the constant $C_T$ there which does not depend on $N$) can grow in $N$ much slower than any power of $N$ and since we already
   explained that $A_N$ decays in $N$ as $N^{-\del}$ for some $\del>0$, Theorem \ref{thm2.6} follows.

\end{document}